\documentclass[journal]{IEEEtran}
\usepackage[nocompress]{cite}
\usepackage{amsmath}
\usepackage{amssymb}
\newtheorem{theorem}{Theorem}

\usepackage{footnote}
\makesavenoteenv{tabular}
\usepackage{algorithm}
\usepackage{algorithmic}

\usepackage{mleftright}
\usepackage{epstopdf}

\usepackage{tikz}
\usetikzlibrary{shapes,arrows,fit,positioning,shadows,calc}
\usetikzlibrary{plotmarks}
\usetikzlibrary{decorations.pathreplacing}
\usetikzlibrary{patterns}
\usetikzlibrary{automata}
\usepackage{pgfplots}
\pgfplotsset{compat=newest}
\usepackage{adjustbox}
\usepackage{caption}
\begin{document}
\title{Compressed Randomized UTV Decompositions for Low-Rank Approximations and Big Data Applications}

\author{\IEEEauthorblockN{Maboud F. Kaloorazi$^{\dagger}$, and %\IEEEauthorrefmark{1},
Rodrigo C. de Lamare$^{\dagger,\ddagger}$} \\%\IEEEauthorrefmark{2}}
\vspace{5mm}
\IEEEauthorblockA{$^\dagger$ Centre for Telecommunications Studies (CETUC)}\\
Pontifical Catholic University of Rio de Janeiro, Brazil \\
$^\ddagger$ Department of Electronics, University of York, United Kingdom \\
E-mail: \texttt{\{kaloorazi,delamare\}@cetuc.puc-rio.br}}

%\vspace{-5cm}

\IEEEtitleabstractindextext{%
\begin{abstract}
Low-rank matrix approximations play a fundamental role in numerical
linear algebra and signal processing applications. This paper
introduces a novel rank-revealing matrix decomposition algorithm
termed Compressed Randomized UTV (CoR-UTV) decomposition along with
a CoR-UTV variant aided by the power method technique. CoR-UTV is
primarily developed to compute an approximation to a low-rank input
matrix by making use of random sampling schemes. Given a large and
dense matrix of size $m\times n$ with numerical rank $k$, where $k
\ll \text{min} \{m,n\}$, CoR-UTV requires a few passes over the
data, and runs in $O(mnk)$ floating-point operations. Furthermore,
CoR-UTV can exploit modern computational platforms and,
consequently, can be optimized for maximum efficiency. CoR-UTV is
simple and accurate, and outperforms reported alternative methods in
terms of efficiency and accuracy. Simulations with synthetic data as
well as real data in image reconstruction and robust principal
component analysis applications support our claims.
\end{abstract}

\begin{IEEEkeywords}
Matrix computations, low-rank approximations, UTV decomposition, randomized algorithms, dimension reduction, matrix decomposition, image reconstruction, robust PCA.
\end{IEEEkeywords}}

\maketitle
\IEEEdisplaynontitleabstractindextext
% \fi
\IEEEpeerreviewmaketitle

\section{Introduction}
\label{sec:intr} \IEEEPARstart{L}{ow-rank} matrix approximations,
that is, approximating a given matrix by one of lower rank, play an
increasingly important role in signal processing and its
applications. Such compact representation which retains most
important information of a high-dimensional matrix can provide a
significant reduction in memory requirements, and more importantly,
computational costs when the latter scales, e.g., according to a
high-degree polynomial, with the dimensionality. Matrices with
low-rank structures have found many applications in background
subtraction \cite{WPMGR2009, BZ2014, KaDeDSP17, RahmaniAtiaHighP17},
system identification \cite{FazelPST13}, IP network anomaly
detection \cite{MMG2013,KaDeICASSP17}, latent variable graphical
modeling, \cite{Chandrasekaran12}, ranking and collaborative
filtering, \cite{Srebro2005}, subspace clustering
\cite{SoltanolkotabiEC2014, RahmaniAtiaCoP17, Oh2017}, adaptive,
sensor and multichannel signal processing
\cite{ifir2005,intadap2005,jio2007,jiomvdr2008,mwfccm2008,DeSa2009,jiols2010},
\cite{ccmjio2010,sjidf2010,jiomimo2011,jiostap2011,barc2011,uwbccm2011,wlmwf2012,wljio2014,dfjio2014,rdrb2015,mserjidf2015,dfalrd2016},
biometrics \cite{VictorBS02, WrightFace09},  statistical process
control and multidimensional fault identification \cite{Jackson91,
Dunia1998}, quantum state tomography \cite{Gross2010}, and DNA
microarray data \cite{Troyanskaya01missingvalue}.

Singular value decomposition (SVD) \cite{GolubVanLoan96} and the
rank-revealing QR (RRQR) decomposition \cite{Chan87, GuEisenstat96}
are among the most commonly used algorithms for computing a low-rank
approximation of a matrix. On the other hand, a UTV decomposition
\cite{StewartURV92, StewartULV93} is a compromise between the SVD
and the RRQR decomposition having the virtues of both: UTV
\textit{i)} is computationally more efficient than the SVD, and
\textit{ii)} provides information on the numerical null space of the
matrix (RRQR does not explicitly furnish the null space information)
\cite{StewartURV92,StewartULV93,Stewart98,Hansen98}. Given a matrix
$\bf A$, the UTV algorithm computes a decomposition ${\bf A
=UTV}^T$, where ${\bf U}$ and ${\bf V}$ have orthonormal columns,
and ${\bf T}$ is triangular (either upper or lower triangular).
These deterministic algorithms, however, are computationally
expensive for large data sets. Furthermore, standard techniques for
their computation are challenging to parallelize in order to utilize
advanced computer architectures \cite{Demmel97, HMT2009, Gu2015}.
Recently developed algorithms for low-rank approximation based on
random sampling schemes, however, have been shown to be remarkably
computationally efficient, highly accurate and robust, and are known
to outperform the traditional algorithms in many practical
situations \cite{FriezeKVS04, DrineasKM06, Sarlos06, Rokhlin09,
HMT2009,Gu2015}. The power of randomized algorithms lies in the
facts that \textit{i)} they are computationally efficient, and
\textit{ii)} their main operations can be optimized for maximum
efficiency on modern computational platforms.

\subsection{Contributions}

Inspired by recent developments, this paper presents a novel
randomized rank-revealing algorithm termed compressed randomized UTV
(CoR-UTV) decomposition \cite{corutv2018}. Given a large and dense
rank-$k$ matrix ${\bf A}$ of size $m \times n$, the CoR-UTV
algorithm computes a low-rank approximation $\hat{\bf A}_\text{CoR}$
of $\bf A$ such that
\begin{equation}
\hat{\bf A}_\text{CoR}={\bf UTV}^T,
\label{eq_contri}
\end{equation}
where ${\bf U}$ and ${\bf V}$ have orthonormal columns, and ${\bf
T}$ is triangular (either upper or lower, whichever is preferred).
CoR-UTV only requires a few passes through data, for a matrix stored
externally, and runs in $O(mnk)$ floating-point operations (flops).
The operations of the algorithm involve matrix-matrix
multiplication, the QR and RRQR decompositions. Due to recently
developed Communication-Avoiding QR algorithms
\cite{DemmGHL12,DemGGX15, DuerschGu2017}, which can perform the
computations with optimal/minimum communication costs, CoR-UTV can
be optimized for peak machine performance on modern architectures.
We provide a theoretical analysis for CoR-UTV, that is, the
rank-revealing property of the algorithm is proved, and upper bounds
on the error of the low-rank approximation are given.

Furthermore, we apply CoR-UTV to treat an image reconstruction
problem, as well as to solve the robust principal component analysis
(robust PCA) problem \cite{WPMGR2009,CSPW2009, CLMW2009}, i.e., to
decompose a given matrix with grossly corrupted entries into a
low-rank matrix plus a sparse matrix of outliers, in applications of
background subtraction in surveillance video, and shadow and
specularity removal from face images.

\subsection{Notation}

Bold-face upper-case letters are used to denote matrices. For a
matrix $\bf A$, ${\|{\bf A}\|_0}$, ${\|{\bf A}\|_1}$, ${\|{\bf
A}\|_2}$, ${\|{\bf A}\|_F}$, and ${\|{\bf A}\|_*}$ denote the
$\ell_0$-norm, the $\ell_1$-norm, the spectral norm, the Frobenius
norm, and the nuclear norm, respectively. $\sigma_j(\bf A)$ and
$\sigma_\text{min}(\bf A)$ denote the $j$-th largest and the
smallest singular value of $\bf A$, respectively. The numerical
range and numerical null space of $\bf A$ are denoted by
$\mathcal{R}({\bf A})$ and $\mathcal{N}({\bf A})$, respectively. The
symbol $\mathbb{E}$ denotes expected value with respect to random
variables, and the dagger $\dagger$ denotes the Moore-Penrose
pseudo-inverse.

The remainder of this paper is structured as follows. In Section
\ref{secRelatW}, we introduce the mathematical model of the data and
discuss related works. In Section \ref{secSOR}, we describe our
proposed method, which also includes a variant that uses the power
iteration scheme in detail.  Section \ref{secAnalysis} presents our
theoretical analysis. In Section \ref{secRobustPCA}, we develop an
algorithm for robust PCA using CoR-UTV. In Section \ref{secNumExp},
we present and discuss our numerical experimental results, and our
concluding remarks are given in Section \ref{secCon}.

\section{Mathematical Model and Related Works}
\label{secRelatW}

Given a matrix ${\bf A} \in \mathbb R^{m \times n}$, where $m \ge
n$, with numerical rank $k$, its singular value decomposition (SVD)
\cite{GolubVanLoan96} is defined as:
\begin{equation}
\begin{aligned}
{\bf A} = & {\bf U}_\text{A}{\bf \Sigma}_\text{A}{\bf V}_\text{A}^T \\ = & \underbrace{\begin{bmatrix} {{\bf U}_k \quad {\bf U}_0} \end{bmatrix}}_{{\bf U}_\text{A} \in \mathbb R^{m \times n}}
  \underbrace{\begin{bmatrix}
       {\bf \Sigma}_k & 0  \\
       0 & {\bf \Sigma}_0
  \end{bmatrix}}_{{\bf \Sigma}_\text{A} \in \mathbb R^{n \times n}}
  \underbrace{\begin{bmatrix}{{\bf V}_k \quad {\bf V}_0} \end{bmatrix}^T}_{{\bf V}_\text{A}^T \in \mathbb R^{n \times n}},
\label{eqSVD}
\end{aligned}
\end{equation}
where ${\bf U}_k \in \mathbb R^{m \times k}$, ${\bf U}_0 \in \mathbb R^{m \times n-k}$ have orthonormal columns, ${\bf \Sigma}_k \in \mathbb R^{k \times k}$ and
${\bf \Sigma}_0 \in \mathbb R^{n-k \times n-k}$ are diagonal matrices containing the singular values, i.e., ${\bf \Sigma}_k=\text{diag}(\sigma_1, ..., \sigma_k)$   and ${\bf \Sigma}_0 =\text{diag}(\sigma_{k+1}, ..., \sigma_n)$, and ${\bf V}_k \in \mathbb R^{n \times k}$ and ${\bf V}_0 \in \mathbb R^{n \times n-k}$ have orthonormal columns. $\bf A$ can be written as ${\bf A} = {\bf A}_k+{\bf A}_0$, where ${\bf A}_k = {\bf U}_k{\bf \Sigma}_k{\bf V}_k^T$, and ${\bf A}_0 = {\bf U}_0{\bf \Sigma}_0{\bf V}_0^T$. The SVD  constructs the optimal rank-$k$ approximation ${\bf A}_k$ to ${\bf A}$, \cite{EckartYoung36, Mirsky60} i.e.,
\begin{equation}
\begin{aligned}
& \underset{\text{rank}({\bf B})\le k}{\text{minimize}}
&& \|{\bf A} - {\bf B}\|_2 = \|{\bf A} - {\bf A}_k\|_2 = \sigma_{k+1}.
\end{aligned}
\label{eq_SVDL2}
\end{equation}
\begin{equation}
\begin{aligned}
&\underset{\text{rank}({\bf B})\le k}{\text{minimize}}
&&\|{\bf A} - {\bf B}\|_F = \|{\bf A} - {\bf A}_k\|_F = \sqrt{\sum_{j=k+1}^{n}{\sigma_j^2}}.
\end{aligned}
\label{eq_SVDF}
\end{equation}

In this paper we focus on the matrix $\bf A$ defined above.

The SVD is highly accurate and yields detailed information on
singular subspaces and singular values. However, it is prohibitive to
compute for large data sets. Moreover, standard techniques for its
computation are challenging to parallelize in order to take advantage of modern
computational environments \cite{Demmel97, HMT2009, Gu2015}.
An economic version of the SVD is the partial SVD based on Krylov subspace methods, such as the Lanczos and Arnoldi algorithms, which constructs an approximate SVD of an input matrix, for instance $\bf A$, at a cost $O(mnk)$. The partial SVD, however, suffers from two drawbacks. First, inherently, it is numerically unstable \cite{CalvettiRS94,GolubVanLoan96, Demmel97}. Second, it does not lend itself to parallel implementations \cite{HMT2009, Gu2015}, which makes it unsuitable for modern computational architectures.

Another widely used algorithm for low-rank approximations considered as a relatively economic alternative to the SVD is the RRQR decomposition \cite{Chan87}. The RRQR is a special QR decomposition with column pivoting (QRCP), which reveals the numerical rank of the input matrix. Given the matrix $\bf A$, it takes the following form:
\begin{equation}
{\bf A}{\bf P} = {\bf Q}{\bf R}= {\bf Q}
  \begin{bmatrix}
       {\bf R}_{11} & {\bf R}_{12}  \\
       {\bf 0} & {\bf R}_{22}
  \end{bmatrix},
\label{equTwolem1}
\end{equation}
where ${\bf P}$ is a permutation matrix, ${\bf Q}\in \mathbb R^{m \times n}$ has orthonormal columns, ${\bf R} \in \mathbb R^{n \times n}$ is upper triangular where ${\bf R}_{11} \in \mathbb R^{k \times k}$ is well-conditioned with
$\sigma_\text{min}({\bf R}_{11})= O(\sigma_k)$, and the $\ell_2$-norm of ${\bf R}_{22}
\in \mathbb R^{n-k \times n-k}$ is sufficiently small, i.e.,
$\|{\bf R}_{22}\|_2= O(\sigma_{k+1})$ (here we have written the reduced QR decomposition, where the silent columns and rows of $\bf Q$ and $\bf R$, respectively, have been removed). If there is an additional requirement that
the $\ell_2$-norm of ${\bf R}_{11}^{-1}{\bf R}_{22}$ is small, i.e., a low order
polynomial in $n$, this decomposition is called ``strong RRQR decomposition" \cite{GuEisenstat96}. The rank-$k$ approximation to $\bf A$ is then computed as follows:
\begin{equation}
\hat{\bf A}_\text{RRQR} = {\bf Q}(:,1:k){\bf R}(1:k,:){\bf P}^T,
\label{A_hat_RRQR}
\end{equation}
where we have used MATLAB notation to indicate submatrices, i.e., ${\bf Q}(:,1:k)$
denotes the first $k$ columns of $\bf Q$, and ${\bf R}(1:k,:)$ denotes the first $k$ rows of $\bf R$.

A UTV decomposition \cite{StewartURV92, StewartULV93} is a compromise between the SVD and QRCP, which has the virtues of both. For the matrix $\bf A$,
it takes the form:
\begin{equation}
{\bf A =UTV}^T
\end{equation}
where ${\bf U}\in \mathbb R^{m \times n}$ and ${\bf V}\in \mathbb R^{n \times n}$ have orthonormal columns, and ${\bf T}$ is triangular. If ${\bf T}$ is upper triangular, the decomposition is called URV decomposition \cite{StewartURV92}:
\begin{equation}
{\bf A}= {\bf U}\begin{bmatrix}
       {\bf T}_{11} & {\bf T}_{12}  \\
       0 & {\bf T}_{22}
       \end{bmatrix}{\bf V}^T.
\label{equ12}
\end{equation}
If ${\bf T}$ is lower triangular, the decomposition is called ULV decomposition \cite{StewartULV93}:
\begin{equation}
{\bf A}= {\bf U}\begin{bmatrix}
       {\bf T}_{11} & 0  \\
       {\bf T}_{21} & {\bf T}_{22}
       \end{bmatrix}{\bf V}^T.
\label{equ13}
\end{equation}

The URV and ULV decompositions are collectively referred to as UTV decompositions \cite{GolubVanLoan96, Hansen98}, and are performed by reduction of the matrix $\bf A$ using unitary transformations to upper and lower triangular forms, respectively.
If there is a well-defined gap in the singular value spectrum of $\bf A$, i.e., $\sigma_k \gg \sigma_{k+1}$, the UTV decompositions are said to be rank-revealing in the sense that the numerical rank $k$ is revealed in the triangular submatrix ${\bf T}_{11} \in \mathbb R^{k \times k}$ \eqref{equ12}, \eqref{equ13}, and the $\ell_2$-norm of off-diagonal submatrices,
$[{\bf T}_{12}^T \quad {\bf T}_{22}^T]^T$ and $[{\bf T}_{21} \quad {\bf T}_{22}]$,
are of the order $\sigma_{k+1}$ \cite{StewartURV92, StewartULV93, FierroHanHan99}, i.e.,
\begin{equation}
\begin{aligned}
&\ \sigma_\text{min}({\bf T}_{11})= O(\sigma_k), \\
&\ \|[{\bf T}_{12}^T \quad {\bf T}_{22}^T ]^T\|_2 = O(\sigma_{k+1}),\\
&\ \|{\bf T}_{21} \quad {\bf T}_{22}\|_2 = O(\sigma_{k+1}).
\end{aligned} %\label{equ15}
\end{equation}

QRCP and UTV decompositions, however, provide highly accurate approximation
to $\bf A$, but they suffer from two drawbacks. First, they are expensive to
compute in terms of arithmetic costs, i.e., $O(mn^2)$. Second, methods for their computation are challenging to parallelize and, as a result, they can not exploit modern architectures \cite{Demmel97, HMT2009, Gu2015}.

Recently developed algorithms for low-rank approximations based on randomization  \cite{FriezeKVS04, DrineasKM06, Sarlos06, Rokhlin09, HMT2009, Gu2015,TrYUC17} have attracted significant attention due to the facts that \textit{i)} they are computationally efficient, and \textit{ii)} they readily lend themselves to a parallel implementation to exploit advanced computational platforms.

The methods in \cite{FriezeKVS04, DrineasKM06,DeshpandeV2006,
RudelsonV07} first sample columns of a given matrix with a
probability proportional to either their magnitudes or leverage
scores, representing the matrix in a compressed form. Further
computations are then performed on the submatrix using deterministic
algorithms such as the SVD and the QR decomposition with column
pivoting \cite{GolubVanLoan96} to obtain the final low-rank
approximation. Sarl{\'{o}}s \cite{Sarlos06} proposed a different
method based on results of the well-known Johnson-Lindenstrauss (JL)
lemma \cite{JL84}. He showed that projecting the data matrix onto a
structured random subspace, i.e., random linear combinations of rows
of the matrix, can render a good approximation to a low-rank matrix.
The works in \cite{NelsonNg2013,ClarWood2017} further advanced
Sarl{\'{o}}s's idea and constructed a low-rank approximation based
on subspace embedding. Rokhlin et al. \cite{Rokhlin09} proposed to
apply a random Gaussian embedding matrix in order to reduce the
dimension of the data matrix. The low-rank approximation was then
obtained through computations using the classical techniques on the
reduced-sized matrix.

Halko et al. \cite{HMT2009} proposed two algorithms based on the
randomized sampling schemes for computing an approximate SVD of an
input matrix. Their first algorithm, called randomized SVD (R-SVD),
projects the matrix onto a low-dimensional subspace using a random
matrix, capturing most attributes of the data. Further computations
are then performed on the reduced-size matrix through the QR
decomposition and the SVD to give the approximation. Gu
\cite{Gu2015} applied a slightly modified version of the R-SVD
algorithm to improve subspace iteration methods, and presents a new
error analysis. The second method proposed by Halko et al.
\cite[Section 5.5]{HMT2009} was a \textit{single-pass} algorithm,
i.e., it required only one pass over the data, to compute a low-rank
approximation. For the matrix ${\bf A}$, the decomposition, which we
call two-sided randomized SVD (TSR-SVD), is computed as described in
Alg. \ref{Alg2}.

\begin{algorithm}
\caption{Two-Sided Randomized SVD (TSR-SVD)}
\renewcommand{\algorithmicrequire}{\textbf{Input:}}
\begin{algorithmic}[1]
\REQUIRE ~~ % ：Input
 Matrix $\ {\bf A} \in \mathbb R^{m \times n}$, integers $k$ and $\ell$.
\renewcommand{\algorithmicrequire}{\textbf{Output:}}
\REQUIRE ~~ A rank-$\ell$ approximation.
  \STATE Draw random matrices ${\bf \Phi}_1 \in \mathbb R^{n \times \ell}$ and
  ${\bf \Phi}_2 \in \mathbb R^{m \times \ell}$;
  \STATE Compute ${\bf Y}_1 = {\bf A}{\bf \Phi}_1$ and ${\bf Y}_2 = {\bf A}^T{\bf \Phi}_2$ in a single pass through $\bf A$;
  \STATE Compute QR decompositions ${\bf Y}_1 = {\bf Q}_1{\bf R}_1$, ${\bf Y}_2 = {\bf Q}_2{\bf R}_2$;
  \STATE Compute ${\bf B}_\text{approx} = {\bf Q}_1^T{\bf Y}_1({\bf Q}_2^T
  {\bf \Phi}_1)^\dagger $;
  \STATE Compute an SVD ${\bf B}_\text{approx} = \widetilde{\bf U}
  \widetilde{\bf \Sigma} \widetilde{\bf V}$;
  \STATE ${\bf A} \approx ({\bf Q}_1 \widetilde{\bf U})\widetilde{\bf \Sigma}
  ({\bf Q}_2 \widetilde{\bf V})^T$.
\end{algorithmic}\label{Alg2}
\end{algorithm}

In Alg. \ref{Alg2}, ${\bf Q}_1 \widetilde{\bf U} \in \mathbb R^{m
\times \ell}$ and ${\bf Q}_2  \widetilde{\bf V} \in \mathbb R^{n
\times \ell}$ are approximations to the left and right singular
subspace of $\bf A$, respectively, $\widetilde{\bf \Sigma} \in
\mathbb R^{\ell \times \ell}$ contains an approximation to the first
$\ell$ singular values of $\bf A$, and ${\bf B}_\text{approx}$ is an
approximation to ${\bf B} = {\bf Q}_1^T{\bf A}{\bf Q}_2$.

TSR-SVD, however, gives a very poor approximation compared to the
optimal SVD due to the single-pass strategy. The reason behind is,
mainly, poor approximate basis drawn from the row space of $\bf A$,
i.e., ${\bf Q}_2$. Furthermore, for a general input matrix, the
authors \cite{HMT2009} do not provide neither theoretical error
bounds nor numerical experiments for the TSR-SVD.

The work in \cite{KaDeDSP17} proposed a rank-revealing decomposition
algorithm based on randomized sampling schemes; the matrix ${\bf A}$
is compressed through pre- and post-multiplication by approximate
orthonormal bases for $\mathcal{R}({\bf A})$ and $\mathcal{R}({\bf
A}^T)$ obtained via randomization, columns of the reduced matrix
and, accordingly, the bases are permuted, the low-rank approximation
is then given by projecting the compressed matrix back to the
original space. The work in \cite{MFKDeTSP18} proposed a randomized
algorithm termed subspace-orbit randomized SVD (SOR-SVD) to compute
a \textit{fixed-rank approximation} of an input matrix. SOR-SVD,
first, alternately projects the matrix onto its column and row
space. Next, orthonormal bases for $\mathcal{R}({\bf A})$ and
$\mathcal{R}({\bf A}^T)$ are approximated. The matrix is then
transformed into a lower dimensional space using the approximate
bases. Finally, a truncated SVD is carried out on the transformed
data, and the low-rank approximation is given by projecting the
small projected data back to the original space.

This work was developed by drawing inspiration from the
rank-revealing algorithm proposed in \cite{KaDeDSP17}, and also
SOR-SVD in \cite{MFKDeTSP18}. Our analysis was inspired by the work
in \cite{MFKDeTSP18}.

\section{Compressed Randomized UTV Decompositions}
\label{secSOR}

In this section, we present a rank-revealing decomposition algorithm
powered by the randomized sampling schemes termed compressed
randomized UTV (CoR-UTV) decomposition, which computes a low-rank
approximation of a given matrix. We focus on the matrix $\bf A$ with
$m \ge n$, where the CoR-UTV algorithm, in the form of
\eqref{eq_contri}, produces a middle matrix $\bf T$, which is upper
triangular, i.e., URV decomposition. The modifications required for
a corresponding CoR-UTV algorithm for the other case, where $m < n$
that produces a lower triangular middle matrix $\bf T$ i.e., ULV
decomposition, is straightforward. For a theoretical comparison of
the URV and ULV decompositions see \cite{Hansen98, FierroHan97} and
the references therein. We also present a version of CoR-UTV with
power iteration, which improves the performance of the algorithm at
an extra computational cost.

Given the matrix ${\bf A}$ and an integer $k\le \ell<
\text{min}\{m,n\}$, the basic version of CoR-UTV is computed as
follows: using a random number generator, we form a matrix ${\bf
\Psi} \in \mathbb R^{n \times \ell}$ with entries drawn independent,
identically distributed (i.i.d.) from the standard Gaussian
distribution. We then compute the matrix product:
\begin{equation}
{\bf C}_1 = {\bf A}{\bf \Psi},
\label{eq_C1_1st}
\end{equation}
where ${\bf C}_1 \in \mathbb R^{m \times \ell}$ is, in fact, a
projection onto the subspace spanned by columns of ${\bf A}$. Having
${\bf C}_1$, we form ${\bf C}_2 \in \mathbb R^{n \times \ell}$:
\begin{equation}
{\bf C}_2 = {\bf A}^T{\bf C}_1,
\label{eqAC1}
\end{equation}
where ${\bf C}_2$ is, in fact, a projection onto the subspace spanned by rows of ${\bf A}$. Using a QR decomposition, we factor ${\bf C}_1$ and ${\bf C}_2$ such that:
\begin{equation}
{\bf C}_1 = {\bf Q}_1{\bf R}_1,  \quad \text{and} \quad {\bf C}_2 =
{\bf Q}_2{\bf R}_2,
\end{equation}
where ${\bf Q}_1$ and ${\bf Q}_2$ are approximate bases for
$\mathcal{R}({\bf A})$ and $\mathcal{R}({\bf A}^T)$, respectively.
We now compress $\bf A$ through left and right multiplications by
the orthonormal bases obtained, forming the matrix ${\bf D} \in
\mathbb R^{\ell \times \ell}$:
\begin{equation}
{\bf D}={\bf Q}_1^T{\bf A}{\bf Q}_2,
\label{eqM}
\end{equation}
We then compute a QRCP of ${\bf D}$:
\begin{equation}
{\bf D} = \widetilde{\bf Q}\widetilde{\bf R}\widetilde{\bf P}^T.
\label{eq_Drrqr}
\end{equation}
Finally, we form the CoR-UTV-based low-rank approximation of $\bf A$:
\begin{equation}
\hat{\bf A}_\text{CoR}= {\bf UTV}^T,
\label{eq_T_basic}
\end{equation}
where ${\bf U}={\bf Q}_1 \widetilde{\bf Q} \in \mathbb R^{m \times
\ell}$ and ${\bf V}={\bf Q}_2  \widetilde{\bf P} \in \mathbb R^{n
\times \ell}$ construct approximations to the $\ell$ leading left
and right singular vectors of $\bf A$, respectively, and ${\bf
T}=\widetilde{\bf R}\in \mathbb R^{\ell \times \ell}$ is upper
triangular with diagonals approximating the first $\ell$ singular
values of $\bf A$. The CoR-UTV algorithm is presented in Alg.
\ref{Alg_first}.

\begin{algorithm}
\caption{Compressed Randomized UTV (CoR-UTV)}
\renewcommand{\algorithmicrequire}{\textbf{Input:}}
\begin{algorithmic}[1]
\REQUIRE ~~ % ：Input
 Matrix $\ {\bf A} \in \mathbb R^{m \times n}$, integers $k$ and $\ell$.
\renewcommand{\algorithmicrequire}{\textbf{Output:}}
\REQUIRE ~~ A rank-$\ell$ approximation.
  \STATE Draw a standard Gaussian matrix ${\bf \Psi}\in \mathbb R^{n\times \ell}$;
  \STATE Compute ${\bf C}_1 = {\bf A}{\bf \Psi}$; \\
  \STATE Compute ${\bf C}_2 = {\bf A}^T{\bf C}_1$;
  \STATE Compute QR decompositions ${\bf C}_1 = {\bf Q}_1{\bf R}_1$, ${\bf C}_2 = {\bf Q}_2{\bf R}_2$;
  \STATE Compute ${\bf D}={\bf Q}_1^T{\bf A}{\bf Q}_2$;
  \STATE Compute the QRCP ${\bf D} = \widetilde{\bf Q}\widetilde{\bf R}\widetilde{\bf P}^T$;
  \STATE Form the CoR-UTV-based low-rank approximation of $\bf A$:
  $\hat{\bf A}_\text{CoR}= {\bf UTV}^T$; ${\bf U}={\bf Q}_1 \widetilde{\bf Q},
  {\bf T}=\widetilde{\bf R}$,${\bf V}={\bf Q}_2\widetilde{\bf P}^T$.
\end{algorithmic}\label{Alg_first}
\end{algorithm}

CoR-UTV requires three passes over the data, for a matrix stored
externally, but it can be altered to revisit the data only once. To
this end, the compressed matrix $\bf D$, equation \eqref{eqM}, can
be computed by making use of available matrices as follows:
\begin{equation}
{\bf D}{\bf Q}_2^T{\bf \Psi} = {\bf Q}_1^T{\bf A}{\bf Q}_2{\bf Q}_2^T{\bf \Psi}.
\end{equation}
where both sides of currently unknown ${\bf D}$ are postmultiplied by ${\bf Q}_2^T{\bf \Psi}$. Having defined ${\bf A}\approx {\bf A}{\bf Q}_2{\bf Q}_2^T$ and ${\bf C}_1 = {\bf A}{\bf \Psi}$, an approximation to ${\bf D}$ is obtained by:
\begin{equation}
{\bf D}_\text{approx} = {\bf Q}_1^T{\bf C}_1({\bf Q}_2^T{\bf \Psi})^\dagger.
\label{eq_approx}
\end{equation}

The key differences between CoR-UTV and TSR-SVD are as follows:
\begin{itemize}
\item CoR-UTV uses a sketch of the input matrix in order to project it onto its row space, i.e., equation \eqref{eqAC1}. This \textit{i)} significantly improves the quality of the approximate basis ${\bf Q}_2$ and, as a result, the quality of the approximate right singular subspace of $\bf A$ compared to that of TSR-SVD, which uses a random matrix for the projection, and \textit{ii)} allows \eqref{eq_approx} to provide a highly accurate approximation to \eqref{eqM}.%\looseness-1

\item CoR-UTV applies a column-pivoted QR decomposition to
$\bf D$, i.e., equation \eqref{eq_Drrqr}, whereas TSR-SVD uses an
SVD to factor the compressed matrix. This, as explained later on, reduces the computational costs of CoR-UTV compared to TSR-SVD.
\end{itemize}

The key difference between CoR-UTV and SOR-SVD \cite{MFKDeTSP18},
however, lies in the computation of the compressed matrix $\bf D$;
SOR-SVD applies a truncated SVD and, as result, gives a rank-$k$
approximation to $\bf A$, while CoR-UTV employs a column-pivoted QR
decomposition and returns a rank-$\ell$ approximation. Nevertheless,
the SVD is computationally more expensive than the column-pivoted
QR, and standard techniques to computing it are challenging to
parallelize \cite{Demmel97, HMT2009, Gu2015}. While recently
developed column-pivoted QR algorithms use randomization, which can
factor a matrix with optimal/minimum communication cost
\cite{DuerschGu2017, MartinssonHQRRP2017}. This can substantially
reduce the computational costs of decomposing the compressed matrix,
compared to the SVD, when it does not fit into fast memory.

CoR-UTV may be sufficiently accurate for matrices whose singular
values display some decay, however in applications where the data
matrix has a slowly decaying singular values, it may produce
singular vectors and singular values that deviate significantly from
the exact ones (computed by the SVD). Thus, we incorporate $q$ steps
of a power iteration \cite{Rokhlin09,HMT2009} to improve the
accuracy of the algorithm in these circumstances. Given the matrix
${\bf A}$, and integers $k\le \ell< n$ and $q$, the resulting
algorithm is described in Alg. \ref{Alg3}.

\begin{algorithm}
\caption{CoR-UTV with Power Method}
\renewcommand{\algorithmicrequire}{\textbf{Input:}}
\begin{algorithmic}[1]
\REQUIRE ~~ % ：Input
 Matrix $\ {\bf A} \in \mathbb R^{m \times n}$,
integers $k$, $\ell$ and $q$.
\renewcommand{\algorithmicrequire}{\textbf{Output:}}
\REQUIRE ~~ A rank-$\ell$ approximation.
  \STATE Draw a standard Gaussian matrix ${\bf C}_2 \in \mathbb R^{n \times \ell}$;
  \FOR{$i=$ 1: $q+1$}
   \STATE Compute ${\bf C}_1 = {\bf A}{\bf C}_2$; \\
   \STATE Compute ${\bf C}_2 = {\bf A}^T{\bf C}_1$;
  \ENDFOR \\
  \STATE Compute QR decompositions ${\bf C}_1 = {\bf Q}_1{\bf R}_1$,
  ${\bf C}_2 = {\bf Q}_2{\bf R}_2$;
  \STATE Compute ${\bf D}={\bf Q}_1^T{\bf A}{\bf Q}_2$ or ${\bf D}_\text{approx} = {\bf Q}_1^T{\bf C}_1({\bf Q}_2^T{\bf C}_2)^
  \dagger $;
  \STATE Compute a QRCP ${\bf D} = \widetilde{\bf Q}\widetilde{\bf R}\widetilde{\bf P}^T$ or ${\bf D}_\text{approx} = \widetilde{\bf Q}\widetilde{\bf R}\widetilde{\bf P}^T$;
  \STATE Form the CoR-UTV-based low-rank approximation of $\bf A$:
  $\hat{\bf A}_\text{CoR}= {\bf UTV}^T$; ${\bf U}={\bf Q}_1 \widetilde{\bf Q},
  {\bf T}=\widetilde{\bf R}$,${\bf V}={\bf Q}_2\widetilde{\bf P}^T$.
\end{algorithmic}\label{Alg3}
\end{algorithm}
Notice that to compute CoR-UTV when the power method is employed, a
non-updated ${\bf C}_2$ must be used to form ${\bf
D}_\text{approx}$.

\section{Analysis of CoR-UTV Decompositions}
\label{secAnalysis}

In this section, we provide an analysis of the performance of
CoR-UTV, the basic version in Alg. \ref{Alg_first} and the one that
uses the power method in Alg. \ref{Alg3}. In particular, we discuss
the rank-revealing property of the algorithm, and provide upper
bounds on the error of the low-rank approximation for CoR-UVT.

We borrow material from \cite{MFKDeTSP18} since the two algorithms,
CoR-UTV and SOR-SVD, have a few steps similar. However, the key
difference is that these randomized algorithms employ different
deterministic decomposition methods to factor the input matrix. We
discuss that CoR-UTV is computationally cheaper and, moreover, can
exploit advanced computer architectures better than SOR-SVD.

\subsection{Rank-Revealing Property}

To prove that CoR-UTV is rank-revealing, it is required to show that
\textit{i}) the $\bf T$ factor of the decomposition reveals the rank
of $\bf A$, and \textit{ii}) the trailing off-diagonal block of $\bf
T$ is small in $\ell_2$-norm. Furthermore, the relation between the
Gaussian random matrix used and the $\bf T$ factor must be
expressed. To be more precise, the quality of the $k$-th
approximated singular value is to be expressed in terms of
properties of the Gaussian matrix.

In CoR-UTV, the $\bf T$ factor is, in fact, the $\bf R$ factor of a
numerically stable deterministic QRCP of $\bf D$, equation
\eqref{eq_Drrqr}, where $\bf D$ is a compressed version of $\bf A$.
We now write \eqref{eq_Drrqr} as:
\begin{equation}
{\bf D}\widetilde{\bf P}= \widetilde{\bf Q}\widetilde{\bf R}=
\widetilde{\bf Q}\begin{bmatrix}
       \widetilde{\bf R}_{11} & \widetilde{\bf R}_{11}  \\
       {\bf 0} & \widetilde{\bf R}_{22}
       \end{bmatrix}.
\label{eq_Drrqr2}
\end{equation}

Thus, it is guaranteed that $\widetilde{\bf R}_{11}\in \mathbb R^{k
\times k}$ reveals the rank of $\bf D$, i.e.,
$\sigma_\text{min}(\widetilde{\bf R}_{11}) \le \sigma_k({\bf D})$,
and $\|\widetilde{\bf R}_{22}\|_2 \le \sigma_{k+1}({\bf D})$.

Next, we need to show how the singular values of $\bf D$ are related to those
of $\bf A$. We establish this relation by stating a theorem from \cite{Thompson72}.

\begin{theorem}(Thompson \cite{Thompson72})
Let the matrix $\bf A$ have an SVD as defined in \eqref{eqSVD}, and ${\bf D} \in \mathbb R^{\ell \times \ell}$ be any submatrix of $\bf A$.
Then for $j = 1,..., \ell,$ we have
\begin{equation}
\sigma_{j+1}\le \sigma_j({\bf D}) \le \sigma_j.
\label{eqThrThomp}
\end{equation}
\end{theorem}

To prove \eqref{eqThrThomp}, it only suffices to allow $\bf D$ be
${\bf D}={\bf Q}_1^T{\bf A}{\bf Q}_2$, where ${\bf Q}_1 \in \mathbb
R^{m \times \ell}$ and ${\bf Q}_2 \in \mathbb R^{n \times \ell}$ are
orthonormal matrices.

Thus, we will have
\begin{equation}
\begin{aligned}
&\sigma_\text{min}(\widetilde{\bf R}_{11}) \le \sigma_k({\bf D}) \le \sigma_k,\\
&\|\widetilde{\bf R}_{22}\|_2 \le \sigma_{k+1}({\bf D}) \le \sigma_{k+1}.
\end{aligned}
\end{equation}

Now, we furnish the relation of $\sigma_k({\bf D})$ and the standard
Gaussian matrix ${\bf \Psi}$. To this end, first, suppose that the
sample size parameter $\ell$ satisfies
\begin{equation}
2\le p+k\le \ell
\end{equation}
where $p$ is called an oversampling parameter \cite{HMT2009,Gu2015}.
Since ${\bf \Psi}$ has interaction with the right singular vectors $\bf V$ of $\bf A$, i.e., equation \eqref{eq_C1_1st}, we have
\begin{equation}
\widetilde{\bf \Psi} = {\bf V}_\text{A}^T{\bf \Psi}=
       [\widetilde{\bf \Psi}_1^T \quad
       \widetilde{\bf \Psi}_2^T]^T
\end{equation}
where $\widetilde{\bf \Psi}_1$ and $\widetilde{\bf \Psi}_2$ have $\ell-p$ and
$n-\ell+p$ rows, respectively. The following theorem, taken from \cite{MFKDeTSP18},
bounds $\sigma_k({\bf D})$.

\begin{theorem}
Suppose that the matrix $\bf A$ has an SVD defined in \eqref{eqSVD},
$2\le p+k\le \ell$, and the matrix $\bf D$ is formed through step 1 to step 5
of Alg. \ref{Alg_first}. Moreover, assume that $\widetilde{\bf \Psi}_1$ is full row
rank, then we have %\looseness-1
\begin{equation}
\sigma_k \ge \sigma_k({\bf D}) \ge \frac{\sigma_k}{\sqrt{1 + {\|{\widetilde{\bf \Psi}_2}\|_2^2} \|{\widetilde{\bf \Psi}_1}^\dagger\|_2^2
\Big(\frac{\sigma_{\ell-p+1}}{\sigma_k}\Big)^4}},
\label{eqThr1_1UTV}
\end{equation}
and when the matrix $\bf D$ is formed through step 1 to step 7 of Alg. \ref{Alg3}, i.e., the power method is used, we have
\begin{equation}
\sigma_k \ge \sigma_k({\bf D}) \ge \frac{\sigma_k}{\sqrt{1 + {\|{\widetilde{\bf \Psi}_2}\|_2^2} \|{\widetilde{\bf \Psi}_1}^\dagger\|_2^2
\Big(\frac{\sigma_{\ell-p+1}}{\sigma_k}\Big)^{4q+4}}}.
\label{eqThr1_2UTV}
\end{equation}
\label{Thr1}
\end{theorem}

Finally, since the random matrix ${\bf \Psi}$ has the standard
Gaussian distribution, the average-case lower bound on the $k$-th
singular value of CoR-UTV is given in the following theorem, taken
from \cite{MFKDeTSP18}.

\begin{theorem}
With the notation of Theorem \ref{Thr1}, and $\gamma_k = \frac{\sigma_{\ell-p+1}}{\sigma_k}$, for Alg. \ref{Alg_first}, we have
\begin{equation}
\mathbb{E}(\sigma_k({\bf D})) \ge \frac{\sigma_k}
{\sqrt{1 + \nu^2\gamma_k^4}},
\end{equation}
and when the power method is used, Alg. \ref{Alg3}, we have
\begin{equation}
\begin{aligned}
\mathbb{E}(\sigma_k({\bf D})) \ge \frac{\sigma_k}
{\sqrt{1 + \nu^2\gamma_k^{4q+4}}},
\end{aligned}
\label{eqThr3_2nd}
\end{equation}
where $\nu = \nu_1\nu_2$, $\nu_1 = \sqrt{n-\ell+p}+\sqrt{\ell}+7$, and $\nu_2 = \frac{4\text{e}\sqrt{\ell}}{p+1}$.
\label{Thr3}
\end{theorem}
This completes the discussion on the rank-revealing property of the
CoR-UTV algorithm.

\subsection{Low-Rank Approximation}

CoR-UTV efficiently constructs an accurate low-rank approximation to
an input matrix $\bf A$. We provide theoretical guarantees on the
accuracy of these approximations in terms of the Frobenius and
spectral norm. To this end, we first state a theorem from
\cite{MFKDeTSP18}.

\begin{theorem}
Let the matrix $\bf A$ have an SVD as defined in \eqref{eqSVD}, and
${\bf Q}_1 \in \mathbb R^{m \times \ell}$ and ${\bf Q}_2 \in \mathbb
R^{n \times \ell}$ be matrices with orthonormal columns constructed
by means of CoR-UTV, where $1 \le k \le \ell$. Let, furthermore,
${\bf D}_k$ be the best rank-$k$ of ${\bf D}={\bf Q}_1^T{\bf A}{\bf
Q}_2$. Then, we have
\begin{equation}
\begin{aligned}
\|{\bf A} - {\bf Q}_1{\bf D}_k{\bf Q}_2^T\|_F \le &
\|{\bf A}_0\|_F +\|{\bf A}_k - {\bf Q}_1{\bf Q}_1^T{\bf A}_k\|_F \\
& + \|{\bf A}_k - {\bf A}_k{\bf Q}_2{\bf Q}_2^T\|_F,
\end{aligned}
\label{equThrFirst2}
\end{equation}
and
\begin{equation}
\begin{aligned}
\|{\bf A} - {\bf Q}_1{\bf D}_k{\bf Q}_2^T\|_2 \le &
\|{\bf A}_0\|_2 +\|{\bf A}_k - {\bf Q}_1{\bf Q}_1^T{\bf A}_k\|_F \\
& + \|{\bf A}_k - {\bf A}_k{\bf Q}_2{\bf Q}_2^T\|_F.
\end{aligned}
\label{equThrFirst3}
\end{equation}
\label{ThrFirst}
\end{theorem}

Now, we rewrite the CoR-UTV-based low-rank approximation
\eqref{eq_T_basic} as follows:
 \begin{equation}
\hat{\bf A}_\text{CoR} = {\bf Q}_1{\bf D}{\bf Q}_2^T.
\label{eq_T}
\end{equation}
This perfectly makes sense since the column-pivoted QR decomposition, which factors
$\bf D$ is a numerically stable deterministic method \cite{GolubVanLoan96}.
Thus, for $\xi=2, F$, it follows that
\begin{equation}
\|{\bf A}-\hat{\bf A}_\text{CoR}\|_\xi \le \|{\bf A}-
{\bf Q}_1{\bf D}_k{\bf Q}_2^T\|_\xi.
\label{equThrFi}
\end{equation}
This relation holds because ${\bf D}_k$ is the rank-$k$ approximation of ${\bf D}$.
\begin{theorem}
With the notation of Theorem \ref{Thr1}, for $\xi =2, F$, we have
\begin{equation}
\begin{aligned}
\|{\bf A} - \hat{\bf A}_\text{CoR}\|_\xi \le &
\|{\bf A}_0\|_\xi +\|{\bf A}_k - {\bf Q}_1{\bf Q}_1^T{\bf A}_k\|_F \\
& + \|{\bf A}_k - {\bf A}_k{\bf Q}_2{\bf Q}_2^T\|_F.
\end{aligned}
\label{equThrFirs}
\end{equation}
\end{theorem}

Having stated the connection between CoR-UTV and SOR-SVD, we now obtain upper bounds for the CoR-UTV-based low-rank approximation error.

\begin{theorem}
Let the matrix $\bf A$ have an SVD as defined in \eqref{eqSVD}, $2\le p+k\le \ell$, and $\hat{\bf A}_\text{CoR}$ is computed through the basic version of CoR-UTV, Alg. \ref{Alg_first}. Furthermore, assume that ${\widetilde{\bf \Psi}_1}$ is full row rank. Then, for $\xi =2, F$, we have
\begin{equation}
\begin{aligned}
\|{\bf A} - \hat{\bf A}_\text{CoR}\|_\xi \le
& \|{{\bf A}_0}\|_\xi +
\sqrt{\frac{\alpha^2{\|{\widetilde{\bf \Psi}_2}\|_2^2} \|{\widetilde{\bf \Psi}_1}^\dagger\|_2^2}{1 + \beta^2{\|{\widetilde{\bf \Psi}_2}\|_2^2} \|{\widetilde{\bf \Psi}_1}^\dagger\|_2^2}} \\ & + \sqrt{\frac{\eta^2{\|{\widetilde{\bf \Psi}_2}\|_2^2} \|{\widetilde{\bf \Psi}_1}^\dagger\|_2^2}{1 + \tau^2{\|{\widetilde{\bf \Psi}_2}\|_2^2} \|{\widetilde{\bf \Psi}_1}^\dagger\|_2^2}},
\end{aligned}
\label{Thr2_1st}
\end{equation}
where $\alpha = \sqrt{k}\frac{\sigma_{\ell-p+1}^2}{\sigma_k}$, $\beta = \frac{\sigma_{\ell-p+1}^2}{\sigma_1 \sigma_k}$, $\eta=\sqrt{k}\sigma_{\ell-p+1}$ and $\tau=\frac{\sigma_{\ell-p+1}}{\sigma_1} $.

When the power iteration is used, Alg. \ref{Alg3}, $\alpha = \sqrt{k}\frac{\sigma_{\ell-p+1}^2}{\sigma_k}\Big(\frac{\sigma_{\ell-p+1}}{\sigma_k} \Big)^{2q}$, $\beta = \frac{\sigma_{\ell-p+1}^2}{\sigma_1\sigma_k}\Big(\frac{\sigma_{\ell-p+1}}{\sigma_k} \Big)^{2q}$, $\eta=\dfrac{\sigma_k}{\sigma_{\ell-p+1}}\alpha$ and  $\tau=\dfrac{1}{\sigma_{\ell-p+1}}\beta$.
\label{Thr2}
\end{theorem}

Theorem \ref{Thr2} implies that, at least, the error bounds for SOR-SVD hold for CoR-UTV. For a detailed error analysis of the SOR-SVD algorithm, see \cite{MFKDeTSP18}.

The random matrix ${\bf \Psi}$ has the standard Gaussian distribution, we thus
present a theorem that establishes average error bounds on the CoR-UTV-based
low-rank approximation.

\begin{theorem}
With the notation of Theorem \ref{Thr2}, and $\gamma_k = \frac{\sigma_{\ell-p+1}}{\sigma_k}$, for the basic version of CoR-UTV, Alg. \ref{Alg_first},
we have
\begin{equation}
\begin{aligned}
\mathbb{E} \|{{\bf A} - \hat{\bf A}_\text{CoR}}\|_\xi \le \mspace{5mu}\|
{{\bf A}_0}\|_\xi + (1 + \gamma_k) \sqrt{k}\nu\sigma_{\ell-p+1},
\end{aligned}
\end{equation}
and when the power method is used, Alg. \ref{Alg3}, we have
\begin{equation}
\begin{aligned}
\mathbb{E} \|{{\bf A} - \hat{\bf A}_\text{CoR}}\|_\xi \le \mspace{5mu}\|{{\bf A}_0}\|_\xi + (1 + \gamma_k) \sqrt{k}\nu\sigma_{\ell-p+1}\gamma_k^{2q},
\end{aligned}
\end{equation}
where $\nu$ is defined in Theorem \ref{Thr3}.
\label{Thr4}
\end{theorem}

This completes the discussion on the low-rank approximation error
bounds for the CoR-UTV algorithm.

\subsection{Computational Complexity}
\label{secComComplex}

The computational cost of any algorithm involves \textit{i)}
arithmetic, i.e., the number of floating-point operations, and
\textit{ii)} communication, i.e., synchronization and data movement
either through levels of a memory hierarchy or between parallel
processing units \cite{DemmGHL12}. On advanced computers, for a
large data matrix which is stored externally, the communication cost
becomes substantially more expensive compared to the arithmetic
\cite{DemmGHL12, Dongarra17}. Therefore, developing new algorithms
or redesigning existing algorithms to solve a problem in hand with
minimum communication cost is highly desirable.

The advantage of algorithms based on randomization over their
classical counterparts lies in the fact that \textit{i)} they
operate on a reduced-size version of the data matrix rather than a
matrix itself, resulting in a reduction of flops, and \textit{ii)}
they can be organized to exploit modern architectures, performing a
decomposition with minimum communication cost.

To factor the matrix $\bf A$, CoR-UTV of Alg. \ref{Alg_first} incurs
the following costs (we only consider high-order terms):
\begin{itemize}
\item Step 1 (generating ${\bf \Psi}$) costs $n\ell$.
\item Step 2 (forming ${\bf C}_1$) costs $2mn\ell$.
\item Step 3 (forming ${\bf C}_2$) costs $2mn\ell$.
\item Step 4 (QR decompositions) costs $2m\ell^2 + 2n\ell^2$.
\item Step 5 (forming ${\bf D}$) costs $2m\ell^2 + 2mn\ell$. If the matrix $\bf D$ is approximated by ${\bf D}_\text{approx}$ of equation \eqref{eq_approx} in this step, the cost would be $2m\ell^2 + 2n\ell^2 +3\ell^3$.
\item Step 6 (performing QRCP) costs $\dfrac{8}{3}\ell^3$.
\item Step 7 (forming the left and right approximate bases) costs $2m\ell^2 + 2n\ell$.
\end{itemize}
Summing up the costs in Steps 1 to 7, we obtain:
\begin{equation}
C_\text{CoR-UTV} \sim 3\ell C_\text{mult}+6m\ell^2 + n\ell(2\ell+3)+ \dfrac{8}{3}\ell^3,
\label{equC1_corutv}
\end{equation}
or
\begin{equation}
C_\text{CoR-UTV} \sim 2\ell C_\text{mult}+6m\ell^2 + n\ell(4\ell+3)+ \dfrac{17}{3}\ell^3,
\label{equC1_corutv_app}
\end{equation}
when the compressed matrix $\bf D$ is approximated by ${\bf
D}_\text{approx}$, where $C_\text{mult}$ is the cost of a
matrix-vector multiplication with $\bf A$ or ${\bf A}^T$. The first
terms of the right-hand sides of \eqref{equC1_corutv} and
\eqref{equC1_corutv_app}, resulting from multiplying $\bf A$ and
${\bf A}^T$ with the corresponding matrices dominate the costs, and
the sample size parameter $\ell$ is typically close to the minimal
rank $k$. When CoR-UTV employs the power method, it requires $2q+3$
passes over the data (for a matrix stored out-of-core) with the
following operation count:
\begin{equation}
C_\text{CoR-UTV} \sim (2q+3)\ell C_\text{mult}+6m\ell^2 + n\ell(2\ell+3)+ \dfrac{8}{3}\ell^3.
\label{equC2_corutv}
\end{equation}
For the case in which the compressed matrix $\bf D$ is approximated
by ${\bf D}_\text{approx}$, CoR-UTV requires $2q+2$ passes over the
data, and the flop count satisfies
\begin{equation}
C_\text{CoR-UTV} \sim (2q+2)\ell C_\text{mult}+6m\ell^2 + n\ell(4\ell+3)+ \dfrac{17}{3}\ell^3.
\label{equC2_corutv_app}
\end{equation}
The CoR-UTV, TSR-SVD and SOR-SVD algorithms except for matrix-matrix
multiplications, which are readily parallelizable perform two QR
decompositions on matrices of size $m\times \ell$ and $n\times
\ell$. Demmel et al. \cite{DemmGHL12} recently developed
communication-avoiding sequential and parallel QR decomposition
algorithms that perform the computations with optimal communication
costs. Hence, this step of all three algorithms can be implemented
efficiently. In addition, CoR-UTV performs one QRCP on an $\ell
\times \ell$ matrix, however TSR-SVD and SOR-SVD perform an SVD on
the $\ell \times \ell$ matrix, which is more expensive than QRCP.
Furthermore, recently developed QRCP algorithms based on
randomization can perform the factorization with minimum
communication costs
\cite{DemGGX15,DuerschGu2017,MartinssonHQRRP2017}, while standard
techniques to compute an SVD are challenging for parallelization
\cite{Demmel97,HMT2009,Gu2015}. As a result, for very large matrices
to be factored on high performance computing architectures, where
the compressed $\ell \times \ell$ matrix does not fit into fast
memory, the execution time to compute CoR-UTV can be substantially
less than those of TSR-SVD and SOR-SVD. This is an advantage of
CoR-UTV over TSR-SVD and SOR-SVD. See \cite{DemmGHL12,DemGGX15} for
a comprehensive discussion on communication costs.

\section{Robust PCA with CoR-UTV}
\label{secRobustPCA} This section describes how to solve the robust
PCA problem using the proposed CoR-UTV method. Principal component
analysis (PCA) \cite{Jackson91} is a widely-used linear
dimensionality reduction technique that tranforms a high-dimensional
data to a low-dimensional subspace which contains most features of
the original data. PCA, however, is known to be very sensitive to
grossly perturbed observations. In order to robustifying PCA against
gross corruption, robust PCA \cite{WPMGR2009,CSPW2009, CLMW2009} was
proposed. Robust PCA represents an input low-rank matrix ${\bf M}
\in \mathbb R^{m \times n}$, whose a fraction of entries being
corrupted, as a linear superposition of a low-rank matrix ${\bf L}$
and a sparse matrix of outliers ${\bf S}$ such as ${\bf M=L+S}$, by
solving the following convex program:
\begin{equation}
\begin{aligned}
&{\text{minimize}_{\bf(L, S)}} \ {\|{\bf L}\|_* + \lambda\|{\bf S}\|_1} \\
&{\text{subject to}} \ {\bf M} = {\bf L} + {\bf S},
\end{aligned}\label{equV1}
\end{equation}
where ${\|\mbox{\bf B}\|_*}  \triangleq \sum_i\sigma_i (\mbox{\bf
B}) $ is the nuclear norm of any matrix $\mbox{\bf B}$,
${\|\mbox{\bf B}\|_1} \triangleq \sum_{ij} |\mbox{\bf B}_{ij}|$ is
the $\ell_{1}$-norm of $\mbox{\bf B}$, and $\lambda>0$ is a tuning
parameter
\cite{rstap2012,sastap2012,locsme2014,saalt2014,als2015,dce2015,damdc2016,okspme2016,bfpeg2011,vfap2012,memd2016},
\cite{spa2008,gbd2013,tds2012,lsmimo2015,mfsic2011,mbthp2014,mbdf2013,mfdf2012,armo2013,sicdma2011,mbsic2011,bfidd2016},
\cite{rmbthp2014,did2014,jpais2012,badstbc2016,wlbd2017,mmimo2013,bbprec2017,1bitidd2018,baplnc2018}.
One efficient method to solve \eqref{equV1} is the method of
augmented Lagrange multipliers (ALM) \cite{Bertsekas1982, YY2009},
which iteratively minimizes the following augmented Lagrangian
function with respect to either variable $\bf L$ or $\bf S$ with the
other one being fixed:
\begin{equation}
\begin{aligned}
\mathcal{L}({\bf L}, {\bf S}, {\bf Y}, \mu) \triangleq %\\&&&
& {\|{\bf L}\|_* + \lambda\|{\bf S}\|_1} + \langle\, {\bf Y}, {\bf M} - {\bf L}\ -{\bf S} \rangle \\
+ & \frac{\mu}{ 2}\|{\bf M - \  L\ -  S}\|_F^2,
\end{aligned} \label{equ4} %\vspace{-0.2cm}
\end{equation}
where $ {\bf Y} \in \mathbb R^{m \times n}$ is the Lagrange
multiplier matrix, and $\mu>0$ is a penalty parameter. The robust
PCA solved by the ALM method is given in Alg. \ref{TableALM-SVD}.

\begin{algorithm}
\caption{Robust PCA solved by ALM}
\renewcommand{\algorithmicrequire}{\textbf{Input:}}
\begin{algorithmic}[1]
\REQUIRE ~~ % ：Input
 Matrix ${\bf M}, \lambda, \mu, {\bf Y}_0 = {\bf S}_0 = 0, j=0$.
\renewcommand{\algorithmicrequire}{\textbf{Output:}}
\REQUIRE ~~ Low-rank plus sparse matrix.
\WHILE {the algorithm does not converge}
        \STATE Compute ${\bf L}_{j+1} = \mathcal{D}_{\mu^{-1}}
        ({\bf M} - {\bf S}_j +\mu^{-1} {\bf Y}_j)$;
        \STATE Compute ${\bf S}_{j+1} = \mathcal{S}_{\lambda\mu^{-1}}
        ({\bf M} - {\bf L}_{j+1} +\mu^{-1} {\bf Y})$;
        \STATE Compute ${\bf Y}_{j+1} = {\bf Y}_j +\mu({\bf M} - {\bf L}_{j+1}
        - {\bf S}_{j+1})$;
\ENDWHILE
\RETURN $\bf L^*$ and $\bf S^*$.
\end{algorithmic}\label{TableALM-SVD}
\end{algorithm}

In Alg. \ref{TableALM-SVD}, for any matrix $\bf B$ with an SVD
defined as ${\bf B} = {\bf U}_\text{B}{\bf \Sigma}_\text{B}{\bf
V}_\text{B}^T$, $\mathcal{D}_\delta ({\bf B})$ refers to a singular
value thresholding operator defined as $\mathcal{D}_\delta ({\bf B})
= {\bf U}_\text{B}\mathcal{S}_\delta ({\bf \Sigma}_\text{B}){\bf
V}_\text{B}^T$, where $\mathcal{S}_\delta (x) =
{\text{sgn}(x)\text{max}}(|x| - \delta, 0)$ is a shrinkage operator
\cite{Hale2008}, and $ \lambda$, $\mu$, ${\bf Y}_0$, and ${\bf S}_0$
are initial values.

The ALM method yields the optimal solution, however its serious
bottleneck is computing a computationally demanding SVD at each
iteration to approximate the low-rank component $\bf L$ of $\bf M$.
To address this concern and to speed up the convergence of the ALM
method, the work in \cite{LLS2011} proposes a few techniques
including predicting the principal singular space dimension, a
continuation technique \cite{Toh2010}, and a truncated SVD by using
PROPACK package \cite{Larsen98}. The modified algorithm
\cite{LLS2011} substantially improves the convergence speed, however
the bottleneck is that the truncated SVD \cite{Larsen98} employed
uses the lanczos algorithm that is inherently unstable and,
moreover, due to the limited data reuse in its operations, it has
very poor performance on modern architectures \cite{CalvettiRS94,
GolubVanLoan96, HMT2009,Gu2015}.

To address this issue, we thus, by retaining the original objective
function proposed in \cite{WPMGR2009,CSPW2009, CLMW2009, LLS2011},
apply CoR-UTV as a surrogate to the truncated SVD to solve the
robust PCA problem. We adopt the continuation technique
\cite{Toh2010, LLS2011}, which increases $\mu$ in each iteration.
The proposed method which is called \texttt{ALM-CoRUTV} is given in
Alg. \ref{TableALM-CoRUTV}.

\begin{algorithm}
\caption{Robust PCA solved by \texttt{ALM-CoRUTV}}
\renewcommand{\algorithmicrequire}{\textbf{Input:}}
\begin{algorithmic}[1]
\REQUIRE ~~ % ：Input
 Matrix ${\bf M}, \lambda, \mu, {\bf Y}_0 = {\bf S}_0 = 0, j=0$.
\renewcommand{\algorithmicrequire}{\textbf{Output:}}
\REQUIRE ~~ Low-rank plus sparse matrix.
\WHILE {the algorithm does not converge}
       \STATE Compute ${\bf L}_{j+1} = \mathcal{C}_{\mu_j^{-1}}
        ({\bf M} - {\bf S}_j +\mu_j^{-1} {\bf Y}_j)$;
        \STATE Compute ${\bf S}_{j+1} = \mathcal{S}_{\lambda\mu_j^{-1}}
        ({\bf M} - {\bf L}_{j+1} +\mu_j^{-1} {\bf Y})$;
        \STATE Compute ${\bf Y}_{j+1} = {\bf Y}_j +\mu_j({\bf M} - {\bf L}_{j+1}
        - {\bf S}_{j+1})$;
        \STATE Update $\mu_{j+1} = \text{max}(\rho\mu_j, {\bar \mu})$;
\ENDWHILE
\RETURN $\bf L^*$ and $\bf S^*$.
\end{algorithmic}\label{TableALM-CoRUTV}
\end{algorithm}

In Alg. \ref{TableALM-CoRUTV}, for any matrix $\bf B$ having a
CoR-UTV decomposition described in Section \ref{secSOR},
$\mathcal{C}_\delta ({\bf B})$ refers to a CoR-UTV thresholding
operator defined as:
\begin{equation}
\mathcal{C}_\delta({\bf B})={\bf U}(:,1:r){\bf T}(1:r,:){\bf V}^T,
\end{equation}
where $r$ is the number of diagonals of $\bf T$ greater than
$\delta$, and $\lambda$, $\mu_0$, ${\bar \mu}$, $\rho$, ${\bf Y}_0$,
and ${\bf S}_0$ are initial values. The main operation of the
\texttt{ALM-CoRUTV} algorithm is computing CoR-UTV in each
iteration, which is efficient in terms of flops, $O(mnk)$, and can
be computed with minimum communication costs; see subsection
\ref{secComComplex}. In subsection \ref{subrpca}, we experimentally
verify that \texttt{ALM-CoRUTV} converges to the exact optimal
solution.

\section{Numerical Experiments}
\label{secNumExp}

In this section, we present the results of some numerical
experiments conducted to evaluate the performance of the CoR-UTV
algorithm for approximating a low-rank input matrix. We show that
CoR-UTV provides highly accurate singular values and low-rank
approximations, and compare our algorithm against several other
algorithms from the literature. We furthermore employ CoR-UTV for
solving the robust PCA problem. The experiments were run in MATLAB
on a desktop PC with a 3 GHz intel Core i5-4430 processor and 8 GB
of memory.

\subsection{Comparison of Rank-Revealing Property \& Singular Values}

We first show that CoR-UTV \textit{i}) is rank revealer, i.e., the
gap in the singular value spectrum of the input matrix is revealed,
and \textit{ii}) provides highly accurate singular values that with
remarkable fidelity track singular values of the matrix. For the
sake of simplicity, we focus on square matrices.

For the randomized algorithms considered, namely CoR-UTV, TSR-SVD,
and SOR-SVD, here and in the next subsection, the results presented
are averaged over 20 trials. Each trial was run with the same input
matrix with an independent draw of the test matrix (or matrices for
TSR-SVD).

We construct two types of matrices of order $10^3$:

\begin{itemize}
\item Matrix 1 (Noisy Low-rank). This matrix is formed by a linear superposition of two matrices ${\bf A} = {\bf A}_1 +{\bf E}$. ${\bf A}_1={\bf U\Sigma V}^T$, where
${\bf U}$ and ${\bf V}$ are random orthonormal matrices, ${\bf \Sigma}$ is a
diagonal matrix containing the singular values $\sigma_i$s that decrease
linearly from $1$ to $10^{-9}$, and $\sigma_{j}=0$ for $j=k+1, ..., 10^3$.
The matrix ${\bf E}$ is a Gaussian matrix normalized to have $\ell_{2}$-norm $\text{gap}\times \sigma_k$. We set the numerical rank $k=20$, and consider two cases:
\begin{itemize}
\item $\text{gap} = 0.01$; \texttt{NoisyLowRank-I}
\item $\text{gap} = 0.1$;  \texttt{NoisyLowRank-II}
\end{itemize}

\item Matrix 2 (Fast Decay). This matrix is formed in a similar way as ${\bf A}_1$ of Matrix 1, but now the diagonals of $\bf \Sigma$ take a different form such that $\sigma_j=1$ for $j=1, ...,k$, and $\sigma_j=(j-k+1)^{-2}$ for $j=k+1, ..., 10^3$ \cite{TrYUC17}. We set the numerical rank $k=10$.
\end{itemize}

We compare the quality of singular values of the matrices considered
computed by our method, described in Section \ref{secSOR}, against
that of alternative rank-revealing methods such as the SVD, QR with
column pivoting (QRCP), UTV, described in Section \ref{secRelatW},
and TSR-SVD (Alg. \ref{Alg2}).

For CoR-UTV and TSR-SVD, we arbitrarily set the sample size
parameter to $\ell=2k$. Both algorithms require the same number of
passes over $\bf A$, either two or $2q+2$ when the power method is
used, to perform a factorization. To compute a UTV decomposition, we
implement the \texttt{lurv} function from \cite{FierroHanHan99}.

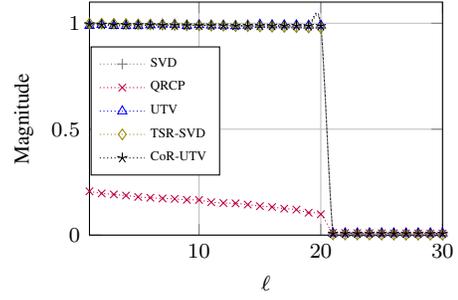
\begin{figure}[t]
\begin{center}
% This file was created by matlab2tikz v0.4.7 running on MATLAB 8.3.
% Copyright (c) 2008--2014, Nico Schlömer <nico.schloemer@gmail.com>
% All rights reserved.
% Minimal pgfplots version: 1.3
% 
% The latest updates can be retrieved from
%   http://www.mathworks.com/matlabcentral/fileexchange/22022-matlab2tikz
% where you can also make suggestions and rate matlab2tikz.
% 
%
% defining custom colors
\usetikzlibrary{positioning,calc}

\definecolor{mycolor1}{rgb}{0.00000,1.00000,1.00000}%
\definecolor{mycolor2}{rgb}{1.00000,0.00000,1.00000}%

\pgfplotsset{every axis label/.append style={font=\footnotesize},
every tick label/.append style={font=\footnotesize}
}

\begin{tikzpicture}[font=\footnotesize] 

\begin{axis}[%
name=ber,
%ymode=log,
width  = 0.53\columnwidth,%5.63489583333333in,
height = 0.35\columnwidth,%4.16838541666667in,
scale only axis,
xmin  = 1,
xmax  = 30,
xlabel= {$\ell$},
xmajorgrids,
ymin=0,
ymax=1.1,
ylabel={Magnitude},
ymajorgrids,
legend entries={SVD, QRCP, UTV, TSR-SVD, CoR-UTV},
legend style={at={(0.37,0.81)},anchor=north east,draw=black,fill=white,legend cell align=left,font=\tiny}
]

\addplot+[smooth,color=gray,densely dotted, every mark/.append style={solid}, mark=+]
table[row sep=crcr]
{
1	0.999859910071192\\
2	0.998926344186047\\
3	0.998389995435944\\
4	0.997035330565599\\
5	0.996019643048321\\
6	0.995498172758966\\
7	0.993965601254841\\
8	0.993227890422083\\
9	0.991956852389459\\
10	0.990816480865243\\
11	0.990102133649130\\
12	0.989178826693723\\
13	0.988248069283534\\
14	0.987066938533269\\
15	0.986311080331270\\
16	0.985050496644538\\
17	0.983926142269320\\
18	0.982703250695982\\
19	0.982082375606026\\
20	0.980686000953933\\
21	0.00970656170546019\\
22	0.00963906629225026\\
23	0.00960779486034902\\
24	0.00958049624186210\\
25	0.00955935020269679\\
26	0.00951676991833767\\
27	0.00949473251032578\\
28	0.00942365827908993\\
29	0.00939544719891800\\
30	0.00938312719531254\\
31	0.00936086246956465\\
32	0.00933392812263468\\
33	0.00932058837403012\\
34	0.00928432357924710\\
35	0.00926313299909757\\
36	0.00923467464064223\\
37	0.00921756672916986\\
38	0.00920703912345985\\
39	0.00916745939393076\\
40	0.00914340956081591\\
}; 
%\addlegendentry{QLP};
\addplot+[smooth,color=purple, densely dotted, every mark/.append style={solid}, mark=x]
  table[row sep=crcr]
  {
1	0.206820574859361\\
2	0.196916084493566\\
3	0.192129454483641\\
4	0.186857525892593\\
5	0.180366206743236\\
6	0.176572348363725\\
7	0.173134730924519\\
8	0.171189163900976\\
9	0.166698909500312\\
10	0.165923612842796\\
11	0.155402140412455\\
12	0.151485882797487\\
13	0.149919607943069\\
14	0.144194537179493\\
15	0.136935477946196\\
16	0.132651359063800\\
17	0.125066143700061\\
18	0.120470269316001\\
19	0.105655798732491\\
20	0.0984396532517840\\
21	0.0106862484657771\\
22	0.00980089662918843\\
23	0.00929471951441145\\
24	0.00902922354336840\\
25	0.00846176166555368\\
26	0.00812726811308465\\
27	0.00795949886015260\\
28	0.00772821574140790\\
29	0.00761172743970554\\
30	0.00737884143594906\\
}; 
\addplot+[smooth,color=blue, densely dotted, every mark/.append style={solid}, mark=triangle]
  table[row sep=crcr]
  {
1	0.988825879318182\\
2	0.993377947403308\\
3	0.987969722283585\\
4	0.987801811928463\\
5	0.989080334042714\\
6	0.991520252236662\\
7	0.994398955566046\\
8	0.993743316132379\\
9	0.988592468412074\\
10	0.992180984293129\\
11	0.990183815063703\\
12	0.989514122149740\\
13	0.987282183182834\\
14	0.989691597601238\\
15	0.993492129149723\\
16	0.991339650555539\\
17	0.991076075348301\\
18	0.991007854507035\\
19	0.989503305769001\\
20	0.990169014488125\\
21	0.00910007864806501\\
22	0.00914097198608600\\
23	0.00884910187467927\\
24	0.00889174734324939\\
25	0.00900999390850653\\
26	0.00874437776010795\\
27	0.00876630686807530\\
28	0.00898601720279713\\
29	0.00891681364204615\\
30	0.00880294921469756\\
31	0.00896526531294066\\
32	0.00884569383658549\\
33	0.00906463750061721\\
34	0.00906026979147614\\
35	0.00874869536633452\\
36	0.00875274783503898\\
37	0.00889735058405034\\
38	0.00863086802599101\\
39	0.00869871968637504\\
40	0.00874692557769848\\
};

\addplot+[smooth,color=olive,densely dotted, every mark/.append style={solid}, mark = diamond]
  table[row sep=crcr]
{
1	0.998897247746082 \\
2	0.997998453487792\\
3	0.997317888948561\\
4	0.996447993244777\\
5	0.995128220738568\\
6	0.994370062658510\\
7	0.993328281402361\\
8	0.992181720900197\\
9	0.990403129514510\\
10	0.989602468138365\\
11	0.988781726625654\\
12	0.988239180943277\\
13	0.986712398560080\\
14	0.985909870428479\\
15	0.985132017403120\\
16	0.983530799174377\\
17	0.982538085716660\\
18	0.981655690624831\\
19	0.979925541599865\\
20	0.978488105757226\\
21	0.00275869396876070\\
22	0.00245970808412507\\
23	0.00244727497157521\\
24	0.00221068510651262\\
25	0.00196374303405508\\
26	0.00170002884685733\\
27	0.00149977906294614\\
28	0.00143820322978009\\
29	0.00138361211393701\\
30	0.00123321246071887\\
31	0.00119295364588483\\
32	0.00112003114676083\\
33	0.000960806594968429\\
34	0.000772124579095462\\
35	0.000697341567365255\\
36	0.000557815566561644\\
37	0.000363044784523014\\
38	0.000293878659127621\\
39	0.000166686750041834\\
40	3.61736631669917e-05\\
};

\addplot+[smooth,color = black ,densely dotted, every mark/.append style={solid}, mark=star]
  table[row sep=crcr]
  {
1	0.995335384289720 \\
2	0.993163055293883\\
3	0.992396696879855\\
4	0.991662727382377\\
5	0.991507972344595\\
6	0.991212331939530\\
7	0.990664123612140\\
8	0.990056378640501\\
9	0.989948011459552\\
10	0.989689808623067\\
11	0.989559244894732\\
12	0.989362817718944\\
13	0.988616213771794\\
14	0.988575852153454\\
15	0.988458549033583\\
16	0.988276311866020\\
17	0.988220455959684\\
18	0.988124678824978\\
19	0.986876329534589\\
20	0.986713099643169\\
21	0.00707629159503820\\
22	0.00701115466517010\\
23	0.00700660060721888\\
24	0.00700250793437634\\
25	0.00697546998430207\\
26	0.00697034816399171\\
27	0.00690292679252935\\
28	0.00684815205278816\\
29	0.00684640151176439\\
30	0.00684136726627775\\
31	0.00682220550111143\\
32	0.00680850878991404\\
33	0.00680224505312939\\
34	0.00677147885378865\\
35	0.00675423489732671\\
36	0.00674300376324053\\
37	0.00668432614706469\\
38	0.00667665145072393\\
39	0.00665666579128886\\
40	0.00645482529707475\\
};
\end{axis}

\end{tikzpicture}%
\captionsetup{justification=centering,font=scriptsize}
\caption{{Comparison of singular values for \texttt{NoisyLowRank-I}
(The power method is not used, $q=0$, for CoR-UTV). CoR-UTV strongly
reveals the gap in the singular values, as do the SVD, UTV and
TSR-SVD, and provides an excellent approximation to singular values.
QRCP weakly reveals the the gap, and considerably underestimates the
leading singular values of the matrix.}}
\label{fig:SV_M1}       % Give a unique label
\end{center}
\end{figure}

\begin{figure}[t]
\begin{center}
\input{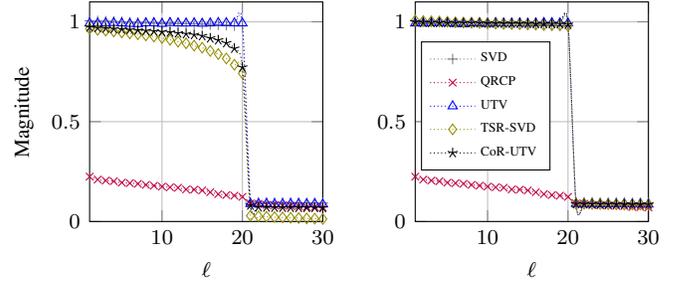}
\captionsetup{justification=centering,font=scriptsize}
\caption{{Comparison of singular values for
\texttt{NoisyLowRank-II}. Left (no power method, $q=0$): CoR-UTV
strongly reveals the gap in the singular values, as do the SVD, UTV
and TSR-SVD, while QRCP only suggests the gap. CoR-UTV outperforms
TSR-SVD in estimating leading and trailing singular values. Right
($q=2$): With two steps of a power method CoR-UTV delivers singular
values as accurate as the optimal SVD. QRCP considerably
underestimates the leading singular values of the matrix.}}
\label{fig:SV_M1_2}       % Give a unique label
\end{center}
\end{figure}
%
%
%\begin{figure}[t]
%\begin{center}
%\input{graph/results_SV_Matrix2}
%\captionsetup{justification=centering,font=scriptsize}
%\caption{\textcolor{blue}{Comparison of singular values for Matrix 2. Left (no power method, $q=0$): Without using the power iterations, CoR-UTV strongly reveals the gap in the singular values, as do the SVD, UTV and TSR-SVD, and gives a very good approximation to singular values. Right ($q=2$): With two steps of a power method CoR-UTV provides singular values as accurate as the optimal SVD. For this matrix, QRCP weakly reveals the the gap, and considerably underestimates the singular values.}}
%\label{fig:SV_M2}       % Give a unique label
%\end{center}
%\end{figure}

The results are shown in Figs. \ref{fig:SV_M1}$-$\ref{fig:SV_M2}. We make the following observations:

\begin{itemize}
\item For all matrices (\texttt{NoisyLowRank-I}, \texttt{NoisyLowRank-II}, Matrix 2), CoR-UTV strongly reveals the numerical rank $k$, as do the SVD, UTV and TSR-SVD, while QRCP weakly reveals the rank of \texttt{NoisyLowRank-I} and Matrix 2, and only suggests the gap in the singular values of \texttt{NoisyLowRank-II}.
\item CoR-UTV, without making use of the power iteration scheme, i.e., $q=0$, provides an excellent approximation to singular values for \texttt{NoisyLowRank-I} and Matrix 2. For \texttt{NoisyLowRank-II}, CoR-UTV outperforms TSR-SVD when $q=0$, in approximating both leading and trailing singular values, and it only requires two steps of the power iteration to deliver singular values as accurate as the optimal SVD. The QRCP algorithm, however, gives a fuzzy approximation to singular values of the input matrices considered.
\end{itemize}

\subsection{Comparison of Low-Rank Approximation}
\subsubsection{Rank-$\ell$ Approximation}
Since CoR-UTV computes a rank-$\ell$ approximation of a given matrix, we first investigate how accurate this approximation is. To this end,
we compute a rank-$\ell$ approximation $\hat{\bf A}_{\text{CoR}}$ for \texttt{NoisyLowRank-I}, \texttt{NoisyLowRank-II}, and Matrix 2 using Alg. \ref{Alg_first} and Alg. \ref{Alg3} for each sample size parameter $\ell$, and calculate the approximation error as:%\looseness-1
\begin{equation}
e_{\ell} = \|{\bf A} - \hat{\bf A}_{\text{CoR}}\|_2.
\label{AppErrorCoR}
\end{equation}
We compare the approximation errors \eqref{AppErrorCoR} against
those produced by the rank-$\ell$ approximation using the SVD, i.e.,
minimal error $\sigma_{\ell+1}$.

Judging from the figures, \textit{i)} when $q=0$, which corresponds to the basic version of CoR-UTV (Alg. \ref{Alg_first}), the approximation is rather poor. \textit{ii)} The error incurred by Alg. \ref{Alg_first} produces an upper bound for the minimal error. \textit{iii)} With only one step of the power iteration $q=1$,
Alg. \ref{Alg3}, the accuracy of the approximation substantially improves,
resulting in an approximation as accurate as the optimal SVD for all three matrices.%\looseness-1
%However, Alg. \ref{Alg_first} provides a highly accurate approximation for Matrix 2 for $\ell \ge 29 $.
%\begin{figure}[t]
%\begin{center}
%\input{graph/results_ErrSVD_Mat1_L2}
%\captionsetup{justification=centering,font=scriptsize}
%\caption{\textcolor{blue}{Comparison of rank-$\ell$ approximation errors of the SVD and CoR-UTV for Matrix 1. When the power method is not utilized ($q=0$), CoR-UTV produces rather poor approximations. When $q=1$ or $q=2$, CoR-UTV approximates the matrices with no loss of accuracy compared to the optimal SVD.}}
%\label{fig:ErrSVD_M1}       % Give a unique label
%\end{center}
%\end{figure}
%
%\begin{figure}[t]
%\begin{center}
%\input{graph/results_ErrSVD_Mat23_L2}
%\captionsetup{justification=centering,font=scriptsize}
%\caption{\textcolor{blue}{Comparison of rank-$\ell$ approximation errors of the SVD and CoR-UTV for Matrix 2. With no power method ($q=0$), CoR-UTV produces rather poor approximation. When $q=1$ or $q=2$, CoR-UTV approximates the matrix with no loss of accuracy compared to the optimal SVD.}}
%\label{ErrSVD_M23}       % Give a unique label
%\end{center}
%\end{figure}

\subsubsection{Rank-$k$ Approximation}
\label{subs_RankKapp}

We now compare the low-rank approximations constructed by our method
against those of the SVD, QRCP, and TSR-SVD. We also include SOR-SVD
\cite{MFKDeTSP18} in our comparison. We have excluded the UTV
algorithm because it has, by far, the worst performance among the
algorithms discussed. This allows us to display the behavior of
other algorithms clearly in the graphs.

To make a fair comparison, we construct a rank-$k$ approximation
${\hat{\bf A}}_\text{out}$ to ${\bf A}$ by each algorithm, and
calculate the error:
\begin{equation}
e_k = \|{\bf A} - \hat{\bf A}_{\text{out}}\|_\xi,
\label{eq_ApprErr}
\end{equation}
where $\xi=F$ for the Frobenius-norm error, and  $\xi=2$ for the
spectral-norm error.

A rank-$k$ approximation for the SVD, QRCP is computed as described
in \eqref{eqSVD} and \eqref{A_hat_RRQR}, respectively. For the
randomized algorithms, TSR-SVD, SOR-SVD, CoR-UTV, however, we
construct a rank-$k$ approximation by varying the sample size
parameter $\ell$, since, as shown, this parameter colors the quality
of approximations. The rank-$k$ approximation by TSR-SVD, Alg.
\ref{Alg2}, is constructed by selecting the first $k$ approximate
singular vectors and corresponding singular values. SOR-SVD
constructs a rank-$k$ approximation of an input matrix, see
\cite{MFKDeTSP18}, and the rank-$k$ approximation by CoR-UTV is
computed as:
\begin{equation}
\hat{\bf A}_{\text{CoR}-k} = {\bf U}(:,1:k){\bf T}(1:k,:){\bf V}^T.
\label{eq_CoR_rank-k}
\end{equation}
For the randomized algorithms, we run the experiment with no power
method ($q=0$), and $q=2$. We make two observations: (1) When $q=0$,
for matrices \texttt{NoisyLowRank-I} and \texttt{NoisyLowRank-II},
as the number of samples increases the performance of CoR-UTV
exceeds that of QRCP, becoming close to optimal performance of the
SVD. For these two matrices, CoR-UTV and SOR-SVD show similar
performances, exceeding the performance of TSR-SVD. For Matrix 2, by
increasing the sample size parameter TSR-SVD and SOR-SVD show
slightly better performance than CoR-UTV, while CoR-UTV outperforms
QRCP. (2) When $q=2$, the errors resulting from CoR-UTV show no loss
of accuracy compared to the optimal SVD. In this case, QRCP has the
poorest performance for all examples.

\subsubsection{Image Reconstruction}

We assess the quality of low-rank approximation by reconstructing a
gray-scale image of a differential gear of size $1280\times 804$,
taken from \cite{DuerschGu2017}, using CoR-UTV, truncated QRCP, and
the truncated SVD by using (widely recommended) PROPACK package
\cite{Larsen98}. The PROPACK function provides an efficient
algorithm to compute a specified number of largest singular values
and corresponding singular vectors of a given matrix by making use
of the Lanczos bidiagonalization algorithm with partial
reorthogonalization, which is suitable for approximating large
low-rank matrices.

The results display the Frobenius-norm approximation error against
the corresponding approximation rank, where the error is calculated
as:
\begin{equation}
e_{\text{approx}} = \|{\bf A} - \hat{\bf A}_{\text{approx}}\|_F,
\label{eq_ApprErr2}
\end{equation}
where $\hat{\bf A}_{\text{approx}}$ is the approximation computed by
each algorithm. %, and Figs. \ref{fig_G20} and \ref{fig_G90} show the
%reconstructed images of the differential gear with $rank=20$ and
%$rank=90$, respectively, using the algorithms mentioned.

For the rank-$20$ approximation, truncated QRCP and CoR-UTV without
power iteration technique produce the poorest reconstruction
qualities. CoR-UTV with one step of power iteration produces a
better result. Truncated SVD and CoR-UTV with two steps of power
iteration, however, appear to have reconstructed images that are
visually identical. For the rank-$90$ approximation, with a careful
scrutiny, fine defects appear in reconstructions by truncated QRCP
and CoR-UTV with $q=0$, while reconstructed images by truncated SVD,
CoR-UTV with $q=1$ and $q=2$ are visually indistinguishable from the
original.

\subsection{Robust Principal Component Analysis}
\label{subrpca}

In this subsection, we experimentally investigate the efficiency and
efficacy of \texttt{ALM-CoRUTV}, described in Table
\ref{TableALM-CoRUTV}, in recovering the low-rank and sparse
components of synthetic and real data. We compare the results
obtained with those of the efficiently implemented inexact ALM
method by \cite{LLS2011}, called \texttt{InexactALM} hereafter.

\subsubsection{Synthetic Matrix Recovery}
\label{subsecSyDataRec}

We form a rank-$k$ matrix $\bf M =L+S$ as a linear combination of a
low-rank matrix ${\bf L} \in \mathbb R^{n \times n}$ and a sparse
error matrix  ${\bf S}\in \mathbb R^{n \times n}$. The matrix ${\bf
L}$ is generated as ${\bf L}={\bf U}{\bf V}^T$, where ${\bf U}$,
${\bf V} \in \mathbb R^{n \times k}$ have standard Gaussian
distributed entries. The error matrix ${\bf S}$ has $s$ non-zero
entries independently drawn from the set $\lbrace$-80, 80$\rbrace$.

We apply the \texttt{ALM-CoRUTV} and \texttt{InexactALM} algorithms
to $\bf M$ to recover ${\bf L}$ and ${\bf S}$. The numerical results
are summarized in Tables \ref{TableCoR1} and \ref{TableCoR2}; Table
\ref{TableCoR1} presents the results where the rank of $\bf L$
$r({\bf L})=0.05\times n$ and $s = \|{\bf S}\|_0=0.05\times n^2$,
and Table \ref{TableCoR2} presents the results for a more
challenging scenario where $r({\bf L})= 0.05\times n$ and $s =
\|{\bf S}\|_0 = 0.10\times n^2$.

In our experiments, we adopt the initial values suggested in
\cite{LLS2011}, and both algorithms are terminated when the
following stopping condition holds:
\begin{equation}
\dfrac{{\|{\bf M}-{\bf L}^{out}-{\bf S}^{out}\|_F}}{{\|{\bf M}\|_F}}< 10^{-5},
\end{equation}
where $({\bf L}^{out}, {\bf S}^{out})$ is the pair of output of
either algorithm. The results of \texttt{ALM-CoRUTV} are reported in
the numerators, and those of \texttt{InexactALM} in the
denominators. In the Tables, $Time(s)$ refers to the computational
time in seconds, $Iter.$ refers to the number of iterations, and
$\zeta$ refers to the relative error defined as ${\|{\bf M}-{\bf
L}^{out}-{\bf S}^{out}\|_F}/{{\|{\bf M}\|_F}}$.

\begin{table}[!htb]
\caption{Numerical results for synthetic matrix recovery for the case
$r({\bf L})=0.05\times n$ and $s = 0.05\times n^2$.}
\begin{tabular}
{p{0.3cm} p{0.3cm} p{0.4cm} p{1.6cm} p{0.4cm} p{0.53cm} p{0.6cm} p{0.3cm} p{0.2cm}}
\noindent\rule{8.4cm}{0.4pt}\\
$n$ & $r(\bf L)$ & $\|{\bf S}\|_0$ & Methods &
$r({\bf L}^*)$ & $\|{\bf S}^*\|_0$ & Time(s) & Iter.& $\zeta$ \\
\noindent\rule{8.4cm}{0.4pt}\\
1000& 50 & 5e4 &
\begin{tabular}{|c p{0.2cm} p{0.7cm} p{0.4cm} p{0.1cm} p{0.9cm}}
\texttt{InexactALM} & 50& 5e4 & 4.1 & 12 & 2.1e-6  \\
\texttt{ALM-CoRUTV} & 50& 5e4 & 0.6 & 12 & 9.6e-6 \\
\end{tabular} \\
\\
2000& 100 & 2e5 &
\begin{tabular}{|c p{0.2cm} p{0.7cm} p{0.4cm} p{0.1cm} p{0.9cm}}
\texttt{InexactALM} & 100& 2e5 & 27.4 & 12 & 2.7e-6  \\
\texttt{ALM-CoRUTV} & 100& 2e5 & 3.7  & 12 & 8.3e-6 \\
\end{tabular} \\
\\
3000& 150 & 45e4 &
\begin{tabular}{|c p{0.2cm} p{0.7cm} p{0.4cm} p{0.1cm} p{0.9cm}}
\texttt{InexactALM} & 150& 45e4 & 75.6 & 12 & 3.1e-6  \\
\texttt{ALM-CoRUTV} & 150& 45e4 & 9.4  & 12 & 8.7e-6 \\
\end{tabular} \\
\\
4000& 200 & 8e5 &
\begin{tabular}{|c p{0.2cm} p{0.7cm} p{0.4cm} p{0.1cm} p{0.9cm}}
\texttt{InexactALM} & 200& 8e5 & 173.3 & 12 & 3.5e-6 \\
\texttt{ALM-CoRUTV} & 200& 8e5 & 20  & 12 & 8.1e-6 \\
\end{tabular} \\
\noindent\rule{8.4cm}{0.4pt}
\end{tabular}
\label{TableCoR1}
\end{table}

\begin{table}[!htb]
\caption{Numerical results for synthetic matrix recovery for the
case $r({\bf L})=0.05\times n$ and $s = 0.1\times n^2$.}
\begin{tabular}
{p{0.3cm} p{0.3cm} p{0.4cm} p{1.6cm} p{0.4cm} p{0.53cm} p{0.6cm} p{0.3cm} p{0.2cm}}
\noindent\rule{8.4cm}{0.4pt}\\
$n$ & $r(\bf L)$ & $\|{\bf S}\|_0$ & Methods &
$r({\bf L}^*)$ & $\|{\bf S}^*\|_0$ & Time(s) & Iter.& $\zeta$ \\
\noindent\rule{8.4cm}{0.4pt}\\
1000& 50 & 1e5 &
\begin{tabular}{|c p{0.2cm} p{0.7cm} p{0.4cm} p{0.1cm} p{0.9cm}}
\texttt{InexactALM} & 50& 1e5 & 4.5 & 13 & 4.4e-6  \\
\texttt{ALM-CoRUTV} & 50& 1e5 & 0.7 & 14 & 9.1e-6 \\
\end{tabular} \\
\\
2000& 100 & 4e5 &
\begin{tabular}{|c p{0.2cm} p{0.7cm} p{0.4cm} p{0.1cm} p{0.9cm}}
\texttt{InexactALM} & 100& 4e5 & 29.2 & 13 & 5.5e-6  \\
\texttt{ALM-CoRUTV} & 100& 4e5 & 4.1  & 14 & 8.9e-6 \\
\end{tabular} \\
\\
3000& 150 & 9e5 &
\begin{tabular}{|c p{0.2cm} p{0.7cm} p{0.4cm} p{0.1cm} p{0.9cm}}
\texttt{InexactALM} & 150& 9e5 & 83.9 & 13 & 6.8e-6  \\
\texttt{ALM-CoRUTV} & 150& 9e5 & 10.9 & 14 & 9.3e-6 \\
\end{tabular} \\
\\
4000& 200 & 16e5 &
\begin{tabular}{|c p{0.2cm} p{0.7cm} p{0.4cm} p{0.1cm} p{0.9cm}}
\texttt{InexactALM} & 200& 16e5 & 189.4 & 13 & 7.8e-6 \\
\texttt{ALM-CoRUTV} & 200& 16e5 & 23.2  & 14 & 9.5e-6 \\
\end{tabular} \\
\noindent\rule{8.4cm}{0.4pt}
\end{tabular}
\label{TableCoR2}
\end{table}

CoR-UTV requires a prespecified rank $\ell$ to perform the factorization. Thus,
we set $\ell=2k$, as a random start, and $q=1$ (one step of a power iteration).
Judging from the results in Tables \ref{TableCoR1} and \ref{TableCoR2}, we make
several observations on \texttt{ALM-CoRUTV}:
\begin{itemize}
\item It successfully detects the exact numerical rank $k$ of the input matrix
in all cases.
\item It provides the exact optimal solution, having the same number of iterations for the first test case, while it requires one more iteration for the second challenging test case, compared to \texttt{InexactALM}.
\item It outperforms \texttt{InexactALM} in terms of runtime, with speedups of up to $8.6$ times.
\end{itemize}

In summary, \texttt{ALM-CoRUTV} exactly recovers the low-rank and
sparse matrices from a grossly corrupted matrix at a much lower cost
compared to \texttt{InexactALM}. However, we expect
\texttt{ALM-CoRUTV} to be faster on multicore and accelerator-based
computers, since CoR-UTV can be computed with minimum communication
cost.

\subsubsection{Background Modeling in Surveillance Video}
\label{subsecBSub}

Extracting the foreground from the background in a video stream is
an increasingly important task in video analysis. This task can be
formulated as a robust PCA problem, where the foreground is
represented by a sparse matrix and the background is represented by
a low-rank matrix.

%\begin{figure}[!htb]
%\centering
%\includegraphics[width=0.48\textwidth,height=4.2cm]{fig/AirporEscalator}
%\caption{Background modeling. Images in columns 1 and 4 are frames of the
%surveillance video of an airport and an escalator, respectively. Images in columns
%2 and 5 are recovered backgrounds ${\bf L}^*$, and columns 3 and 6 correspond
%to foregrounds ${\bf S}^*$ by \texttt{ALM-CoRUTV}.}
%\label{fig_AirEsc}
%\end{figure}

%\begin{figure}[!htb]
%\centering
%\includegraphics[width=0.48\textwidth,height=3.8cm]{fig/HighTram}
%\caption{\textcolor{blue}{Background modeling. Images in columns 1 and 4 are frames of the
%surveillance video of a highway and a tram stop, respectively. Images in columns
%2 and 5 are recovered backgrounds ${\bf L}^*$, and columns 3 and 6 correspond
%to foregrounds ${\bf S}^*$ by \texttt{ALM-CoRUTV}.}}
%\label{fig_HighTram}
%\end{figure}

Here, we apply \texttt{ALM-CoRUTV} to four different surveillance
videos. The first two videos are from \cite{LHGT2004}, and the other
two are from \cite{CDNET14}. The first video consists of 200
grayscale frames of size ${176 \times 144}$, taken in a hall of an
airport. The frames are stacked as columns of a matrix $\bf M$,
forming ${\bf M} \in \mathbb R^{25344 \times 200}$. This video has a
relatively static background. The second video consists of 200
grayscale frames of size ${130 \times 160}$, taken from an escalator
at an airport. We form a matrix ${\bf M} \in \mathbb R^{20800 \times
200}$ by stacking individual frames as its columns. The background
of this video changes due to the moving escalator. The third video
has 200 grayscale frames of size ${240 \times 320}$, taken from a
highway.  We thus form ${\bf M} \in \mathbb R^{76800 \times 200}$.
The fourth video consists of 200 grayscale frames of size ${288
\times 432}$, taken in a tram stop. Therefore, ${\bf M} \in \mathbb
R^{124416 \times 200}$.

In order for CoR-UTV, used in \texttt{ALM-CoRUTV}, to approximate
the low-rank component of real data, we determine the prespecified
rank $\ell$ by making use of the following bound that relates the
numerical rank $k$ of any matrix $\bf B$ with the nuclear and
Frobenius norms \cite{GolubVanLoan96}:
\begin{equation}
\|{\bf B}\|_* \le \sqrt{k}\|{\bf B}\|_F.
\label{equRank}
\end{equation}
We set $\ell=k + p$, where $k$ is the minimum value satisfying \eqref{equRank}, and
$p = 2$ is an oversampling parameter. Again, we set $q=1$ for CoR-UTV.

%Some frames of the surveillance videos with recovered foregrounds
%and backgrounds are displayed in Figs. \ref{fig_AirEsc} and
%\ref{fig_HighTram}. We only show the results of \texttt{ALM-CoRUTV}
%since those produced by \texttt{InexactALM} are visually identical.
%As can be seen the proposed \texttt{ALM-CoRUTV} can successfully
%recover the low-rank and sparse components of the videos. Table
%\ref{TableSix} presents the numerical results. In all four examples,
%\texttt{ALM-CoRUTV} outperforms \texttt{InexactALM} in terms of
%runtime.

\begin{table}[!htb]
\caption{Numerical results for real matrix recovery.}
\begin{tabular}{p{2cm} p{2.5cm} p{2.5cm}}
\noindent\rule{8.1cm}{0.4pt}\\
Dataset &
 \begin{tabular}{p{0.6cm} p{0.5cm} p{0.4cm}}
 \multicolumn{3}{c}{\texttt{InexactALM}} \\
\hline
 Time(s) & Iter. & $\zeta$ \\
 \end{tabular}
 & \begin{tabular}{p{0.6cm} p{0.5cm} p{0.4cm}}
   \multicolumn{3}{c}{\texttt{ALM-CoRUTV}} \\
   \hline
   Time(s) & Iter. & $\zeta$ \\
   \end{tabular}\\
\noindent\rule{8.1cm}{0.4pt}\\
Airport hall $25344 \times 200$ &
\begin{tabular}{p{0.65cm} ccp{0.85cm}c c}
15.4 & 28 & 7.4e-6 & $\mspace{2mu}$ 5.1 & $\mspace{-18mu}$ 28 & 7.1e-6\\
\end{tabular} \\
Escalator $20800 \times 200$ &
\begin{tabular}{p{0.65cm} ccp{0.85cm}c c}
11.9 & 28 & 8.5e-6 & $\mspace{2mu}$ 4.2 & $\mspace{-18mu}$ 28 & 6.2e-6\\
\end{tabular} \\
Highway $76800 \times 200$ &
\begin{tabular}{p{0.65cm} ccp{0.85cm}c c}
53.6 & 28 & 6.5e-6 & $\mspace{2mu}$ 16.2 & $\mspace{-18mu}$ 28 & 8.6e-6\\
\end{tabular} \\
Tram stop $124416 \times 200$ &
\begin{tabular}{p{0.65cm} ccp{0.85cm}c c}
83.6 & 27 & 8.9e-6 & $\mspace{2mu}$ 25.2 & $\mspace{-18mu}$ 28 & 6.3e-6\\
\end{tabular} \\
Yale B01 $32256 \times 64$ &
\begin{tabular}{p{0.65cm} ccp{0.85cm}c c}
4.2 & 26 &7.6e-6& $\mspace{3mu}$ 2.1 & $\mspace{-18mu}$ 26 & 8.6e-6\\
\end{tabular} \\
Yale B02 $32256 \times 64$ &
\begin{tabular}{p{0.65cm} ccp{0.85cm}c c}
4.2 & 26 &6.5e-6& $\mspace{3mu}$ 2.1 & $\mspace{-18mu}$ 26 &  8.3e-6\\
\end{tabular} \\
\noindent\rule{8.1cm}{0.4pt}\\
\end{tabular}
\label{TableSix}
\end{table}

%\begin{figure}[!htb]
%\centering
%\includegraphics[width=0.42\textwidth,height=4cm]{fig/ImYaleCoR}
%\caption{Removing shadows and specularities from face images. Images in columns 1 and 4 are cropped and aligned images of a face under different illuminations. Images in columns 2 and 5 are recovered images by \texttt{ALM-CoRUTV}, and images in columns 3 and 6 are sparse errors corresponding to the removed shadows and specularities.}
%\label{figFaceCoR}
%\end{figure}

\subsubsection{Shadow and Specularity Removal From Face Images}
\label{subsecShRemoval}

Another task in computer vision that fits nicely into the robust PCA
model is removing shadows and specularities from face images; images
of the same face taken under varying illumination can be modeled as
a superposition of a low-rank and sparse components.

In this experiment, we use face images taken from the Yale B face
database \cite{Georghiades2001}.% Each image has the size ${192
%\times 168}$ with a total of 64 different illuminations. The
%individual images are stacked as columns of a matrix ${\bf M} \in
%\mathbb R^{32256 \times 64}$. The recovered images are displayed in
%Fig. \ref{figFaceCoR}. From this figure, we observe that the shadows
%and specularities have been effectively extracted in the sparse
%components by \texttt{ALM-CoRUTV}.
Table \ref{TableSix} summarizes
the numerical results.

We conclude that \texttt{ALM-CoRUTV} can successfully recover the face images
under different illuminations from the dataset studied two times faster than \texttt{InexactALM}.

\section{Conclusion}
\label{secCon}

In this paper, we have presented CoR-UTV, a rank-revealing algorithm
based on the randomized sampling paradigm, for computing a low-rank
approximation of an input matrix. We have presented theoretical
analysis for CoR-UTV, and have shown through numerical experiments
on two classes of matrices that CoR-UTV reveals the numerical rank
better than QRCP, and provides results as good as those of the
optimal SVD. CoR-UTV outperforms QRCP in low-rank approximation, and
when the power method is incorporated, provides results as accurate
as those of the SVD. CoR-UTV is more efficient than the
deterministic SVD, QRCP, UTV, and competing randomized TSR-SVD and
SOR-SVD in terms of arithmetic cost and, moreover, can exploit
advanced computational platforms better by exposing higher levels of
parallelism than all algorithms mentioned. We also applied CoR-UTV
to solve the robust PCA problem via the ALM method. Our studies
demonstrate that the resulting \texttt{ALM-CoRUTV} provides the
exact optimal solution and, moreover, is substantially faster than
efficiently implemented \texttt{InexactALM}.

%\appendices
%\section{The spectral-norm error comparison for rank-$k$ approximation}
%
%This appendix provides numerical results on the performance of the
%algorithms discussed in subsection \ref{subs_RankKapp} for the
%rank-$k$ approximation of Matrix 1 and Matrix 2 in terms of the
%spectral-norm error.

%\begin{figure}[t]
%\begin{center}
%\input{graph/res_LowRank_L2NMat1_I}
%\caption{Comparison of the spectral-norm rank-$k$ approximation error for  \texttt{NoisyLowRank-I}. Left: No power method, $q=0$. Right: $q=2$.}
%\label{fig:LR_L2M1_I}       % Give a unique label
%\end{center}
%\end{figure}
%
%\begin{figure}[t]
%\begin{center}
%\input{graph/res_LowRank_L2NMat1_II}
%\caption{Comparison of the spectral-norm rank-$k$ approximation error for  \texttt{NoisyLowRank-II}. Left: No power method, $q=0$. Right: $q=2$.}
%\label{fig:LR_L2M1_II}       % Give a unique label
%\end{center}
%\end{figure}
%
%\begin{figure}[t]
%\begin{center}
%\input{graph/res_LowRank_L2NMat2}
%\caption{Comparison of the spectral-norm rank-$k$ approximation error for Matrix 2. Left: No power method, $q=0$. Right: $q=2$.}
%\label{fig:LR_L2M2}       % Give a unique label
%\end{center}
%\end{figure}

\bibliographystyle{IEEEtran}
\bibliography{mybibfile}
\end{document}